\RequirePackage{amsmath}
\documentclass[graybox]{svmult}
\usepackage{mathptmx}
\usepackage{helvet}  
\usepackage{courier} 
\usepackage{type1cm} 
\usepackage{makeidx}  
\usepackage{graphicx} 
\usepackage{multicol} 
\usepackage[bottom]{footmisc}
\usepackage{amsmath}
\usepackage{amsfonts}
\usepackage{algorithmic}
\usepackage{algorithm}

\newcommand{\N}{\mathrm{N}}
\newcommand{\NN}{\mathcal{N}}

\newcommand{\statdom}{\Omega_\mathrm{s}}
\newcommand{\rotdom}{\Omega_\mathrm{r}}
\newcommand{\interface}{\Gamma_\mathrm{I}}
\newcommand{\boundrot}{\Gamma_\mathrm{r}}
\newcommand{\boundstat}{\Gamma_\mathrm{s}}

\begin{document}
\title*{Model Order Reduction for Rotating Electrical Machines}
\author{Zeger Bontinck, Oliver Lass, Sebastian Sch\"ops and Oliver Rain}
\institute{Zeger Bontinck \at Technische Universit\"at Darmstadt, Graduate School of Computational Engineering, Dolivostr. 15, 64293 Darmstadt, Germany, \email{bontinck@gsc.tu-darmstadt.de}
\and Oliver Lass \at  Technische Universit\"at Darmstadt, Department of Mathematics, Chair of Nonlinear Optimization, Dolivostr. 15, 64293 Darmstadt, Germany, \email{lass@mathematik.tu-darmstadt.de}
\and Sebastian Sch\"ops \at Technische Universit\"at Darmstadt, Graduate School of Computational Engineering, Dolivostr. 15, 64293 Darmstadt, Germany, \email{schoeps@gsc.tu-darmstadt.de}
\and Oliver Rain \at Robert Bosch GmbH, 70049 Stuttgart, Germany, \email{oliver.rain@de.bosch.com}}
\maketitle
\abstract{The simulation of electric rotating machines is both computationally
expensive and memory intensive. To overcome these costs, model order reduction
techniques can be applied. The focus of this contribution is especially on
machines that contain non-symmetric components. These are usually introduced
during the mass production process and are modeled by small perturbations in
the geometry (e.g., eccentricity) or the material parameters. While model order
reduction for symmetric machines is clear and does not need special treatment,
the non-symmetric setting adds additional challenges. An adaptive strategy
based on proper orthogonal decomposition is developed to overcome these
difficulties. Equipped with an a posteriori error estimator the obtained
solution is certified. Numerical examples are presented to demonstrate the
effectiveness of the proposed method.} \section{Introduction}

\label{sec:intro}
Model order reduction for partial differential equations is a very active field
in applied mathematics. When performing simulations in $2$D or $3$D using the
finite element method (FEM) one arrives at large scale systems that have to be
solved. Projection based model order reduction methods have shown to
significantly reduce the computational complexity when applied carefully. While
being applied to many different fields in physics, the application to rotating
electrical machines is more recent \cite{HC14,LU16,MHCG16,SSSI16}. We will
focus especially on the setting of non-symmetric machines. While the perfect
machine is symmetric and simulations are usually carried out exploiting these
properties, in real life the symmetry is often lost. This is due to
perturbation in the geometry (e.g., eccentricity) and material properties and
requires that the whole machine is simulated and not only a small portion
(e.g., one pole). Hence this leads to an increase in the computational cost.
The aim is to develop an adaptive strategy that is able to collect the required
information systematically. Ideally, the algorithm is able to detect symmetries
and exploits them if present. The greedy algorithm introduced in the context of
the reduced basis method is a possible candidate \cite{PR06,QMN16}. A commonly
used method in engineering and applied mathematics is the method of snapshots,
or proper orthogonal decomposition (POD) \cite{Cha00,HLBR12,Sir87}. We opt for
a combination of the two methodologies. The goal is a fast and efficient
variant that avoids the necessity of an online-offline decomposition. This
allows a broader application since no expensive offline costs have to be
invested. Hence the method will already pay off after one simulation and not
only in the many query context.

Additionally, the developed strategy has to be able to handle the motion of the
rotor. While there are a number of methods to treat the rotation
\cite{DRL85,PC95,SLDP06}, we will assume a constant rotational speed which
allows us to utilize the locked step method \cite{PRS88}. Hence we can avoid
the remeshing which is computationally prohibitive. Moreover, the application
of other approaches should be straight forward.

Efficient simulation tools are a key ingredient when performing optimization or
uncertainty quantification. The combination of model order reduction and
optimization has caught a lot of attention \cite{DH15,GV17,NGMQ13,QGVW17,ZF15}.
Especially in the many query context, where models have to be evaluated
repeatedly and the need for accurate, fast and reliable reduced order models is
high. While we will not look into the application of the reduced order models,
we will develop a strategy that fulfills these needs. By using existing
simulation tools in the adaptive procedure it is possible to insert the
developed strategy into an existing framework and utilize the benefits of the
reduced order model, as shown in \cite{LU16}.

The article is structured as follows: In Section~\ref{sec:prob} the permanent
magnet synchronous machine (PMSM) is introduced. The discretization and the
proposed model order reduction strategy are introduced in
Section~\ref{sec:mor}. Then in Section~\ref{sec:num} numerical experiments are
presented. Lastly, a conclusion is drawn in the last section.

\section{Problem Description}
\label{sec:prob}

\begin{figure}[b]
\centering
 \includegraphics[width=58mm]{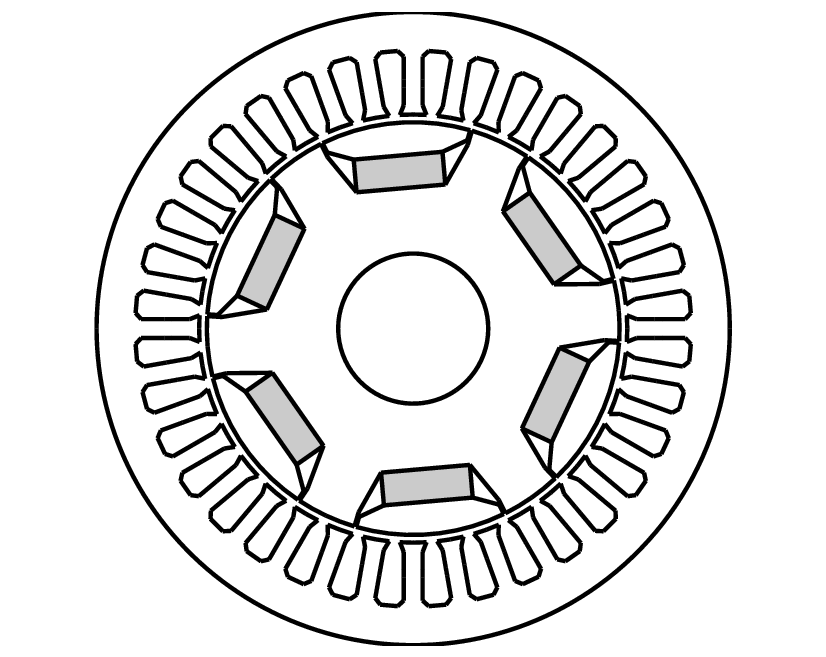}
 \includegraphics[width=58mm]{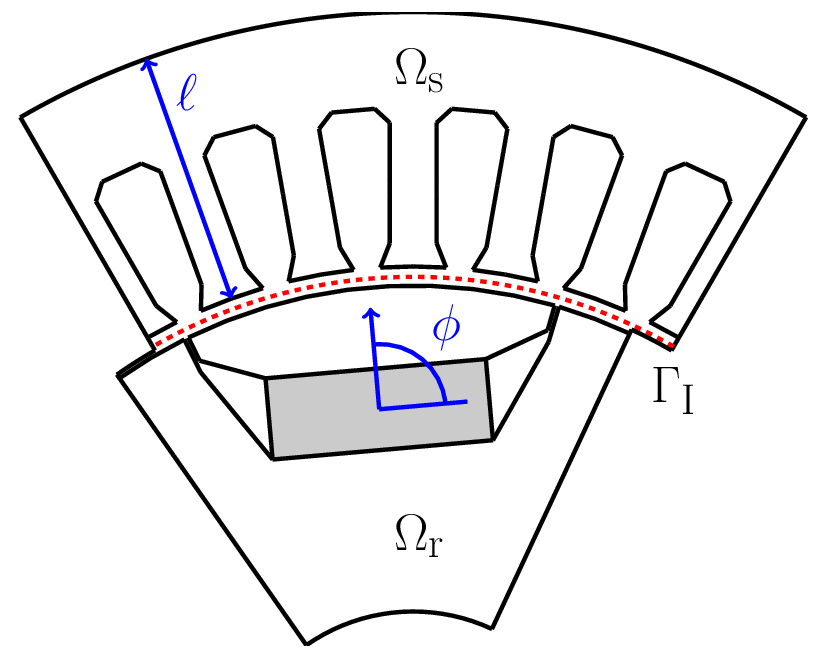}
\caption{Cross sectional view of the full PMSM (left) and detailed view on one pole right.}\label{fig:geom}
\end{figure}

The PMSM under investigation has six slots per pole and a double layered
winding scheme with two slots per pole per phase. The geometry of the full
machine is shown in Figure~\ref{fig:geom} (left) together with a detailed view
on one pole (right). In each pole there is one buried permanent magnet
indicated in gray. The machine has depth $\ell_z = 10\,\mathrm{mm}$. It is
operated at $50\,\mathrm{Hz}$, resulting in a rotational speed of $1000$
revolutions per minute (RPM). The machine is composed of laminated steel with a
relative permeability $\mu_r=500$.
In the following $\statdom$ and $\rotdom$ refer to the stator and rotor domain,
respectively. The whole domain is then referred to as $\Omega = \statdom \cup
\rotdom$ and will be used for simpler notation when appropriate. Additionally,
let us define the interface $\interface=\partial\statdom \cap \partial\rotdom$
(dashed line) in the airgap between the rotor and the stator. Furthermore, we
introduce the boundaries $\boundstat=\partial\statdom\setminus\interface$ and
$\boundrot=\partial\rotdom\setminus\interface$ of the stator and rotor,
correspondingly. In the simulation we will account for the movement of the
rotor, hence we introduce the angle $\vartheta$ that describes the position of
the rotor with respect to the stator. For clarity we will append $\vartheta$ to
the components related to the rotor.

To calculate magnetic vector potential of the machine, the magnetostatic
approximation of the Maxwell's equations has to be solved for both domains.
This implies that the eddy and displacement currents are neglected and one
obtains the semi-elliptic partial differential equations
\begin{equation}\label{eq:mqscont}
    \vec{\nabla}\times\left(\nu\vec{\nabla}\times\vec{A}(\vartheta)\right) =\vec{J}_{\rm src}(\vartheta)-\nabla\times\left(\nu\mathbf{B}_{\rm rem}\right), \quad \mbox{on } \statdom
\end{equation}
with Dirichlet boundary conditions $\vec{A}\times\mathbf n =0$ on $\boundstat$,
with $\mathbf n$ the outer unit normal. The reluctivity is depicted by a scalar
$\nu$ since only linear isotropic materials are considered and nonlinearity is
disregarded since a linearization at a working point is assumed.
$\vec{A}(\vartheta)$ is the magnetic vector potential, $\vec{J}_{\rm
src}(\vartheta)$ represents the imposed source current density, which is
related to the applied currents in the coils, and $\mathbf{B}_{\rm rem}$ the
remanence of the permanent magnets. The applied current density is aligned with
the $z$-direction, whereas the remanence is in the $xy$-plane. It is generally
accepted that machines are adequately modeled in $2$D, meaning that the
magnetic field has no $z$-component: $\vec{B}=(B_x, B_y, 0)$. Since
$\vec{B}=\nabla\times\vec{A}$, one can write $\vec{A} = (0, 0, A_z)$. Hence we
end up with the linear elliptic equation
\begin{equation}
 \label{eq:mqs}
 -\nabla\cdot(\nu\nabla A_z(\vartheta)) = j_{\mathrm{src}}(\vartheta) - j_{\mathrm{pm}}, \quad \mbox{on } \Omega(\vartheta),
\end{equation}
where $j_{\mathrm{src}}$ and $j_{\mathrm{pm}}$ are the $z$-component of $\vec{J}_{\rm src}$ and  $\nabla\times\left(\nu\mathbf{B}_{\rm rem}\right)$, respectively. 
Further, the boundary conditions are given as previously introduced by $A_z = 0$ on $\boundstat$.

Let us next have a look at possible imperfections in the presented geometry and
model. These are introduced during mass production of the PMSMs and are given
by small deviations in the geometry or material properties. While there are
many possible properties that can occur, we focus on two types. On the one hand
we look at imperfection in the material of the permanent magnet, more precisely
we consider deviations in the magnetic field angle $\phi$ of the permanent
magnet \cite{OH13}. The second imperfection we consider is the length $\ell$ of
the teeth in the stator \cite{Cle13,OMNCDH15}. Both quantities are depicted in
Figure~\ref{fig:geom} (right). Perturbations in these quantities may lead to
underperformance of the PMSM.

\section{Model order reduction}
\label{sec:mor}

In this section we discuss a method to accelerate the simulation of the PMSM.
We will start by introducing the finite element discretization. Further, the
realization of the rotation is outlined for the discrete setting. The resulting
linear system of equations are of large scale and hence expensive to solve. In
a next step we will present a model order reduction method based on proper
orthogonal decomposition to speedup the simulations.

\subsection{Finite element discretization}
\label{subsec:fem}

We obtain the discrete version of (\ref{eq:mqs}) by utilizing the finite
element method (FEM). Discretizing $A_z$ by linear edge shape functions
$w_i(x,y)$ one makes the Ansatz
\begin{equation*}
A^\N_z=\sum_{i=1}^\N a_i w_i(x,y)=\sum_{i=1}^\N a_i\frac{\varphi_i(x,y)}{\ell_z}e_z,
\end{equation*}
where $\varphi_i(x,y)$ depicts the nodal finite elements which are associated
with the triangulation of the machine's cross-section and $e_z$ is the unit
vector in $z$-direction. Using the Ritz-Galerkin approach, the $\N$-dimensional
linear discrete system
\begin{equation}
\label{eq:discr}
\mathbf{K}_\nu(\vartheta) \mathbf{a}(\vartheta) =\mathbf{f}(\vartheta)
\end{equation}
is obtained, where $\mathbf{K}_\nu$ are the finite element system matrices, $\mathbf{a}$ depict the degrees of freedom (DoFs) and $\mathbf{f} = \mathbf{j}_{\mathrm{src}} - \mathbf{j}_{\mathrm{pm}}$ the discretized versions of the current densities and permanent magnets. 

To take the motion of the rotor into account in the simulation, we utilize the
locked step method \cite{PRS88,SLDP06}. For the implementation, a contour in
the airgap is defined ($\interface$) which splits the full domain into two
parts: a fixed outer domain connected to the stator $\statdom$ and an inner
domain connected to the rotor $\rotdom$, where the mesh will be rotated. At the
contour the nodes are distributed equidistantly. The angular increments
$\Delta\vartheta$ are chosen so that the mesh of the stator and rotor will
match on the interface. The nodes on the interface are then reconnected leading
to the mesh for the next computation. Using this technique, the rotation angle
$\vartheta$ is given by $\vartheta^k = k\Delta\vartheta$ with $k\in
\mathbb{N}_0$. We can hence partition the discrete unknown $\mathbf a$ in into
a static part, a rotating part and the interface, with dimensions $\N_s$,
$\N_r$ and $\N_I$, respectively.
This idea is a particularization of non-overlapping domain decomposition \cite{AHH11,TW05}.
The linear system (\ref{eq:discr}) can then be written as
\begin{equation}
\label{eq:discr_rot}
 \begin{bmatrix}
   \mathbf K^{ss}_\nu & 0  & \mathbf K^{sI}_\nu \\
   0   & \mathbf K^{rr}_\nu& \mathbf K^{rI}_\nu(\vartheta) \\
   (\mathbf K^{sI}_\nu)^\top & (\mathbf K^{rI}_\nu)^\top(\vartheta) & \mathbf K^{II}_\nu(\vartheta)
 \end{bmatrix}
 \begin{bmatrix}
   \mathbf a_s \\
   \mathbf a_r \\
   \mathbf a_I
 \end{bmatrix}
 =
 \begin{bmatrix}
   \mathbf f^s \\
   \mathbf f^r \\
   \mathbf f^I(\vartheta)
 \end{bmatrix},
\end{equation}
where $\mathbf K^{ss}_\nu$, $\mathbf K^{rr}_\nu$, $\mathbf f^s$ and $\mathbf
f^r$ are the stiffness matrices and right hand sides on the static and moving
part and do no longer depend on the angle $\vartheta$. For the points on the
interface there are two cases. The interface of the static part is independent
of the angle $\vartheta$ and hence we get the corresponding stiffness matrix
$\mathbf K^{sI}_\nu$. For the rotor side we have to perform the shift, this is
indicated with $\vartheta$ in the corresponding stiffness matrix $\mathbf
K^{rI}_\nu$. On the interface also a shift has to be performed hence also here
the corresponding stiffness matrix $\mathbf K^{II}_\nu$ and right hand side
$\mathbf f^I$ are dependent on $\vartheta$. Let us note that it is not required
to reassemble matrices. All of these shifts can be performed by index shift and
hence allow a very efficient implementation. The size of the system does not
change, i.e., we have $\N = \N_s + \N_r + \N_I$.

For completeness let us remark on the domain deformation introduced by $\ell$.
Since we are in the setting of a parametrized shape transformation, i.e., one
$\ell$ for each tooth in the stator we can utilize an affine geometry
preconditioning to describe the transformation. This allows to carry out the
computation on a reference domain and hence avoids again expensive remeshing
and its undesired consequences as e.g., mesh noise. The stiffness matrices and
right hand sides can be written in the form
$$
\mathbf{K}_\nu(\vartheta) = \sum_{j=1}^L \Theta(\ell) \mathbf{K}_\nu^j
\quad\mbox{and}\quad
\mathbf{f}(\vartheta) = \sum_{j=1}^L \Theta(\ell) \mathbf{f}^j,
$$
where $\mathbf{K}_\nu^j$ and $\mathbf{f}^j$ are local matrices. 
For a detailed description including computational procedures, we refer the reader to \cite{RHP08}.

\subsection{Proper orthogonal decomposition}
\label{subsec:pod}

The goal is to generate a reduced order model to accelerate the simulation of
(\ref{eq:discr_rot}). The simulation of the rotation is computationally
expensive since the discretization of (\ref{eq:mqs}) leads to very large linear
systems that need to be solved for every angular position. While in symmetric
machines this can be avoided, in the case of non-symmetric machines a whole
revolution has to be simulated. Hence we require an efficient strategy to
overcome this problem. For this we investigate an adaptive strategy that builds
a surrogate model while performing the simulation and switches to it when the
required accuracy is reached. We want to use information collected over the
rotational angle $\vartheta$ and generate a projection based reduced order
model. In the past, model order reduction methods based on proper orthogonal
decomposition (POD) \cite{Cha00,HLBR12,KV01}, balanced truncation
\cite{AHH11,HSS08} and the reduced basis method \cite{PR06, QMN16, RHP08} have
been developed to speedup the computation. More recently, the POD method has
been successfully applied to rotating machines \cite{HC14,LU16,MHCG16,SSSI16}.
In this work we consider a combination of POD and the reduced basis method. We
will not pursue an online-offline decomposition but rather see the reduced
order model as an accelerator for the simulation.

Let us start by recalling the POD method so we can develop an extension
suitable for the application presented. Let the solution to (\ref{eq:mqs}) be
given in the discrete form, i.e., let $\mathbf a(\vartheta) \in \mathbb R^{\N}$
be the solution to (\ref{eq:discr}) for a fixed angle $\vartheta$. The
snapshots are then given by $\mathbb R^{\N} \ni \mathbf a^k\approx \mathbf
a(k\Delta\vartheta)$ for $k\in\mathcal K$, where $\mathcal K$ is an index set
with elements in $\mathbb N_0$ for which (\ref{eq:discr}) is solved. A POD
basis $\{ \psi_i\}_{i=1}^{\ell}$ is then computed from these snapshots by
solving the following optimization problem:
\begin{equation}
\label{eq:pod}
\left\{\begin{array}{l}
\displaystyle\min_{\psi_1,\ldots,\psi_{\NN} \in  \mathbb R^{\N}}
\sum_{k\in\mathcal K} \,  
\Big|\mathbf a^{k}(\mu) - \sum_{i=1}^{\NN}
{\langle \mathbf a^{k}(\mu),\psi_i\rangle}_\mathrm W\,
\psi_i\Big|_\mathrm W^2\vspace{2mm}\\
\hspace{6.5mm}\text{s.t. } 
{\langle \psi_i,\psi_j \rangle}_\mathrm W = \delta_{ij}
\text{ for } 1 \le i,j \le \NN,
\end{array}\right.
\end{equation}
where $\langle\cdot\,,\cdot\rangle_\mathrm W$ stands for the weighted inner
product in $ \mathbb R^{\N}$ with a positive definite, symmetric matrix
$\mathrm W\in \mathbb R^{\N\times \N}$. The goal is to minimize the mean
projection error of the given snapshots projected onto the subspace spanned by
the POD basis $\psi_i$. By introducing the matrix $\mathbf A_{\mathcal K}$ as
the collection of the snapshots $\mathbf a^k$ with $k\in\mathcal K$ we can
write the operator $\mathbf R$ arising from the optimization problem
(\ref{eq:pod}) as
\begin{equation*}
\mathbf R\psi = \sum_{k\in\mathcal K}\langle \mathbf a^{k}(\mu),\psi\rangle_{\mathrm W}\mathbf a^{k} = \mathbf A_{\mathcal K} \big(\mathbf A_{\mathcal K}\big)^\top \mathrm W \psi
\quad\text{for }\psi \in \mathbb R^{\N}.
\end{equation*}
Then the unique solution to \eqref{eq:pod} is given by the eigenvectors corresponding to the $\NN$ largest eigenvalues of $\mathbf R$, i.e., $\mathbf R \psi_i = \lambda_i \psi_i$ with $\lambda_i > 0$ \cite{GV17}. The operator $\mathbf R$ is largesince it is of dimension $\N$ which we want to reduce. 
Hence it might be better in many cases to set up and solve the eigenvalue problem
\begin{equation*}
 \mathbf A_{\mathcal K}^\top\mathrm W\mathbf A_{\mathcal K} v_i = \lambda_i v_i,\quad i = 1,\ldots,\NN
\end{equation*}
and obtain the POD basis by $\psi_i = 1/\sqrt{\lambda_i}\mathbf A_{\mathcal
K}v_i$. Note that both approach are equivalent and are related by the singular
value decomposition (SVD) of the matrix $\mathrm W^{1/2}\mathbf A_{\mathcal
K}$. While the latter is computationally more efficient, the singular value
decomposition is numerically more stable. A comparison of the different
computations was carried out in \cite{LV13}. For completeness let us state the
POD approximation error given by
\begin{equation}
 \label{eq:pod_approximation_error}
 \sum_{k\in\mathcal K} \,  
\Big|\mathbf a^{k} - \sum_{i=1}^{\NN}
{\langle \mathbf a^{k},\psi_i\rangle}_\mathrm W\,
\psi_i\Big|_\mathrm W^2 = \sum_{i=\NN+1}^d\lambda_i,
\end{equation}
where $d$ is the rank of $\mathbf A_{\mathcal K}$. 
For easier notation we collect the POD basis $\psi_i$ in the matrix $\Psi = [\psi_1,\ldots,\psi_\NN] \in \mathbb R^{\N \times \NN}$. 

After introducing the computation of the POD basis we will now outline the
adaptive approach utilized in this work. It is crucial to minimize the number
of solves involving the FEM discretization in order to obtain a speedup of the
computation. The goal is to push most of the computations in the simulation
onto the reduced order models. However, we have to guarantee that the reduced
order models are accurate in order to obtain reliable results. We will present
a strategy that does not require precomputation as for example in the reduced
basis method but performs the model order reduction during the simulation.
First let us outline how the POD basis is applied to \eqref{eq:discr_rot}. In
the second step we give the details on how to obtain the basis efficiently.

We generate for each part of the machine an individual POD basis. Hence we have
one basis for the stator and one basis for the rotor. The interface between the
stator and the rotor is not reduced, i.e., we work with the FEM ansatz space on
the interface. This is motivated by the observation that the decay of the
eigenvalues on $\interface$ is very slow, which would result in a large POD
basis. Since the FEM space for the interface is usually of moderate dimension,
the gain of using POD would be negligible. We compute the POD basis as solution
to \eqref{eq:pod} utilizing the snapshots $\mathbf a_s^k$ and $\mathbf a_r^k$
to obtain $\Psi^s$ and $\Psi^r$, respectively. 
We then make the ansatz
\begin{equation*}
 \mathbf a^\NN_s = \sum_{i=1}^{\NN_s} \psi_i^s \bar{\mathbf a}_{s,i} = \Psi^s \bar{\mathbf a}_s
 \quad\text{and}\quad
 \mathbf a^\NN_r = \sum_{i=1}^{\NN_r} \psi_i^r \bar{\mathbf a}_{r,i} = \Psi^r \bar{\mathbf a}_r,
\end{equation*}
where the POD coefficients are indicated with a bar. 
By projecting \eqref{eq:discr_rot} onto the subspace spanned by the POD basis, we obtain the reduced order model 
\begin{equation*}
\begin{bmatrix}
  (\Psi^s)^\top\mathbf K^{ss}_\nu\Psi^s & 0  & (\Psi^s)^\top\mathbf K^{sI}_\nu \\
  0   & (\Psi^r)^\top\mathbf K^{rr}_\nu\Psi^r& (\Psi^r)^\top\mathbf K^{rI}_\nu(\vartheta) \\
  (\mathbf K^{sI}_\nu)^\top\Psi^s & (\mathbf K^{rI}_\nu)^\top(\vartheta)\Psi^r & \mathbf K^{II}_\nu(\vartheta)
\end{bmatrix}
\begin{bmatrix}
  \bar{\mathbf a}_s \\
  \bar{\mathbf a}_r \\
  \mathbf a_I
\end{bmatrix}
 = 
\begin{bmatrix}
  (\Psi^s)^\top\mathbf f^s \\
  (\Psi^r)^\top\mathbf f^r \\
  \mathbf f^I(\vartheta)
\end{bmatrix}.
\end{equation*}
In short notation the system will be written as 
\begin{equation}
 \label{eq:discr_rom}
 \mathbf K^\NN_\nu(\vartheta) \bar{\mathbf a} = \mathbf f^\NN(\vartheta)
\end{equation}
with $\bar{\mathbf a}$ the vector of POD coefficients. 
This system is of dimension $\NN_s + \NN_r + \N_I$ and of much smaller dimension as the original system \eqref{eq:discr_rot} which was of dimension $\N$. 

Next we introduce the strategy on how to determine the POD basis adaptively.
The goal is to reduce the computational cost with respect to the rotation. A
full revolution requires $\N_I$ solves of the system \eqref{eq:discr_rot},
i.e., for all $\vartheta_k$ with $k\in \mathcal K := \{0,1,\ldots,\N_I-1\}$. In
the symmetric case it is not required to solve a full rotation but only for
angles that cover one pole, for our particular example this means one sixth,
i.e., $\N_I/6$ solutions are needed. Note, we assume that $\N_I$ is divisible
by $\N_p$, where $\N_p$ is the number of poles of the machine. In the
non-symmetric case this is not possible. The idea is to generate a sequence of
disjoint sets $\mathcal K_i$ with $i = 0,\ldots,\N_{\mathcal K}$. Note that the
sets can also be chosen arbitrarily as long as they fulfill the subset
property. This is required to be able to reuse the already computed snapshots
and hence minimize computational overhead. Ideally the sets $\mathcal K_i$ are
not too large but are large enough to cover the most important features. The
strategy is then as follows: We start by choosing $\mathcal K_0$ and evaluate
\eqref{eq:discr_rot} for $\vartheta_k$ and $k \in \mathcal K_0$. From the
computed solutions/snapshots $A_{\mathcal K} = [\mathbf a^k]_{k\in\mathcal
K_0}$ a POD basis is computed. Then an error estimator $\Delta\mathbf
a(\vartheta_k)$, $k = 1,\ldots,N_I$, is evaluated to determine the maximum
error. If the error is larger than a given tolerance the index $k \in
\{1,\ldots,\N_I\}$ is determined where the maximum error is attained. Determine
the set $\mathcal K_i$ that contains the index $k$ and the snapshot set is
enlarged by adding the new solutions corresponding to $\mathcal K_i$ to the old
ones, i.e., $A_{\mathcal K} = [A_{\mathcal K},[\mathbf a^k]_{\mathcal K_i}]$.
This procedure is repeated until the error $\Delta\mathbf a(\vartheta_k)$, $k =
1,\ldots,N_I$, is below the desired threshold. For stability reasons the sets
$\mathcal K_i$ that have been already used are removed from the list. In this
case the index for the second largest error is used. During the numerical
experiments this scenario never occurred and hence will not be investigated
further. In Algorithm~\ref{alg:pod} the strategy is summarized. This sampling
of the sets $\mathcal K_i$ is similar to the greedy algorithm from the reduced
basis method. The decision to add more than one solution $\mathbf a^k$ at a
time to the snapshot set is to minimize the overhead of evaluating the error
estimator and generating the reduced order models. Since we do not introduce an
online-offline decomposition the computations of error estimator and the
generation of the reduced order models are included in the computational costs.
In our numerical tests the proposed strategy converges very fast, in at most
$6$ iterations. Depending on the asymmetry different sets $\mathcal K_i$ may
preform better. This will be illustrated with numerical experiments. For the
symmetric case we observe that usually only one set is visited (if chosen
appropriately). In other words: the symmetry is correctly detected by the
algorithm. In Section~\ref{sec:num} the specific choices of the sets $\mathcal
K_i$ are described and a comparison of different strategies is performed. The
dimension $\NN$ of the computation of the POD basis is chosen such that
\begin{equation*}
\frac{\sum_{i=1}^{\NN}\lambda_i}{\sum_{i=1}^{d}\lambda_i}\le\varepsilon_{rel}
\end{equation*} 
holds for the stator and rotor independently. This is a popular choice, where a
typical value is $\varepsilon_{rel} = 0.9999$. Note that the denominator can be
computed by $\mathrm{trace}(\mathbf A_{\mathcal K}^\top\mathrm W\mathbf
A_{\mathcal K})$ and hence not all $d$ eigenvalues have to be computed.

\begin{algorithm}
\caption{Adaptive POD}
\label{alg:pod}
\begin{algorithmic}[1]
\REQUIRE{$\mathcal K_{i=0}^{\N_{\mathcal K}}$ and $\varepsilon$ (tolerance) }
\STATE{Choose first set, e.g., $i=0$ and set the snapshot set $\mathbf A_{\mathcal K} = []$}
\STATE{\label{alg:pod:goto} Solve $\mathbf a^k$ for $k \in\mathcal K_i$ and add to $\mathbf A_{\mathcal K}$}
\STATE{Compute POD basis using \eqref{eq:pod}}
\STATE{Evaluate error estimator $\Delta^{\mathrm{rel}}_{\mathbf a}(\vartheta_k)$ for $k = 1,\ldots,\N_I$}
\IF{$\max_k \Delta^{\mathrm{rel}}_{\mathbf a}(\vartheta_k) > \varepsilon$}
\STATE{Determine index $i$ of set $\mathcal K_i$ containing $k$}
\STATE{GOTO \ref{alg:pod:goto}}
\ELSE
\RETURN POD basis and reduced solution $\mathbf a^\ell$
\ENDIF
\end{algorithmic}
\end{algorithm}

Let us now shortly have a look at the error estimator. For this let us recall some basic quantities. 
We define the discrete coercivity constant by
\begin{equation*}
 \alpha(\vartheta) = \inf_{v\in \mathbb R^\N \setminus \{0\}} \frac{\mathbf v^\top \mathbf K_\nu(\vartheta) \mathbf v}{\mathbf v^\top \mathbf W \mathbf v}.
\end{equation*}
Hence the the coercivity constant is given by the smallest eigenvalue such that 
\begin{equation*}
 \mathbf K_\nu(\vartheta) \mathbf v = \lambda \mathbf W \mathbf v
\end{equation*}
is satisfied for $(\lambda,\mathbf v) \in \mathbb R_{+} \times \mathbb R^\N$ and $\mathbf v \neq 0$ \cite{QMN16}.
Further, we define the residual $r(\mathbf a^\NN;\vartheta) = \mathbf f(\vartheta) - \mathbf K_\nu(\vartheta) \mathbf a^\NN$. 
Then the error introduced by the reduced order model in the variable $\mathbf a$ can be characterized by
\begin{equation}
\label{eq:error_estimator_state}
 \|\mathbf a - \mathbf a^\NN\|_{\mathbf W}  \le \Delta_{\mathbf a}(\vartheta) := \frac{\|r(\mathbf a^\NN;\vartheta)\|_{{\mathbf W}^{-1}}}{\alpha(\vartheta)}.
\end{equation}
Additionally, we look at the relative error which might be more interesting in many applications. 
The corresponding error estimator reads
\begin{equation}
\label{eq:error_estimator_state_rel}
 \frac{\|\mathbf a - \mathbf a^\NN\|_{\mathbf W}}{\|\mathbf a^\NN\|_{\mathbf W}}  \le \Delta^{\mathrm{rel}}_{\mathbf a}(\vartheta) := 2\frac{\|r(\mathbf a^\NN;\vartheta)\|_{{\mathbf W}^{-1}}}{\alpha(\vartheta)\|\mathbf a^\NN\|_{\mathbf W}}.
\end{equation}
This is a standard result and can be found in \cite{PR06,RHP08}. In the
numerical realization we set the weight matrix $\mathbf W$ to the $\N$
dimensional identity matrix. Other choices (e.g., $\mathbf W = \mathbf K_\nu$)
are possible but not investigated at this point. Note that the rotation as
introduced in this work does not influence the coercivity constant and hence
the dependence can be omitted. This can also be seen in the numerics, where
only small deviations can be observed which are in the order of discretization.
Hence a very efficient realization is possible, since only one eigenvalue
problem has to be solved. Let us remark that the error is measured with respect
to the finite element solution. It is assumed that the finite element solution
is accurate enough to approximate the solution of the continuous problem. Let
us remark that the computation of the residual norms can be performed very
efficiently using the introduced affine decomposition \cite{DH15}.

\section{Numerical results}
\label{sec:num}

We will now present different numerical results. For this let us specify the
settings. The geometry (Figure~\ref{fig:geom}) is discretized using a
triangular mesh with $56297$ nodes and the interface $\interface$ is
discretized by $900$ equidistant points. Hence one revolution requires the
linear system (\ref{eq:discr}) to be solved $900$ times. Here is where model
order reduction will come into play and speedup the simulation significantly.
All computation are performed on a standard desktop PC using Matlab R2016b.
Throughout our numerical tests we will consider three settings:
\begin{description}
\item[\emph{sym}\hspace{5mm}] Symmetric machine.
\item[\emph{rot}\hspace{6.5mm}] Perturbation of $\phi$ in one permanent magnet by $5^\circ$.
\item[\emph{stat}\hspace{5.25mm}] Perturbation of $\ell$ in one tooth by $0.3$mm.
\item[\emph{rot\_stat}] Perturbation in both $\phi$ and $\ell$. 
\end{description} 

\begin{figure}
    \begin{tabular}{cc}
          \includegraphics[width=58mm]{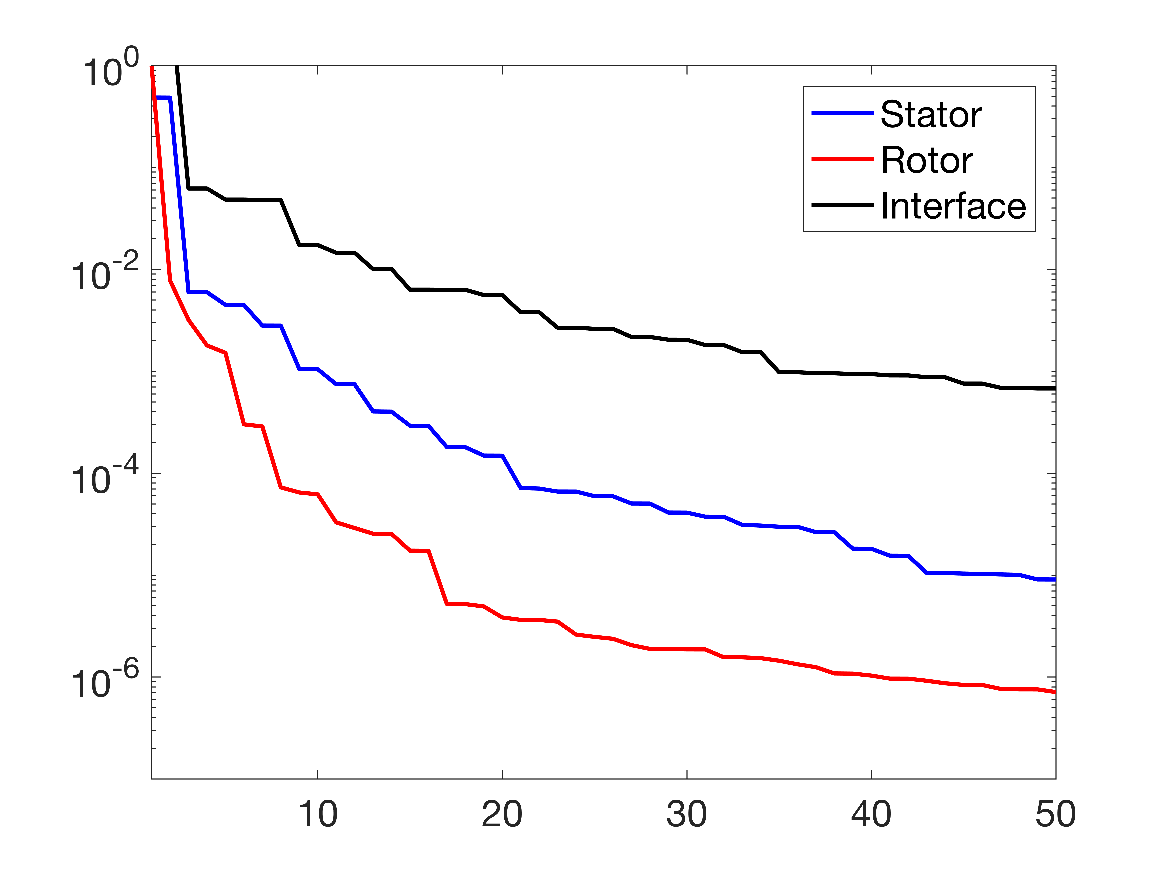}
        & \includegraphics[width=58mm]{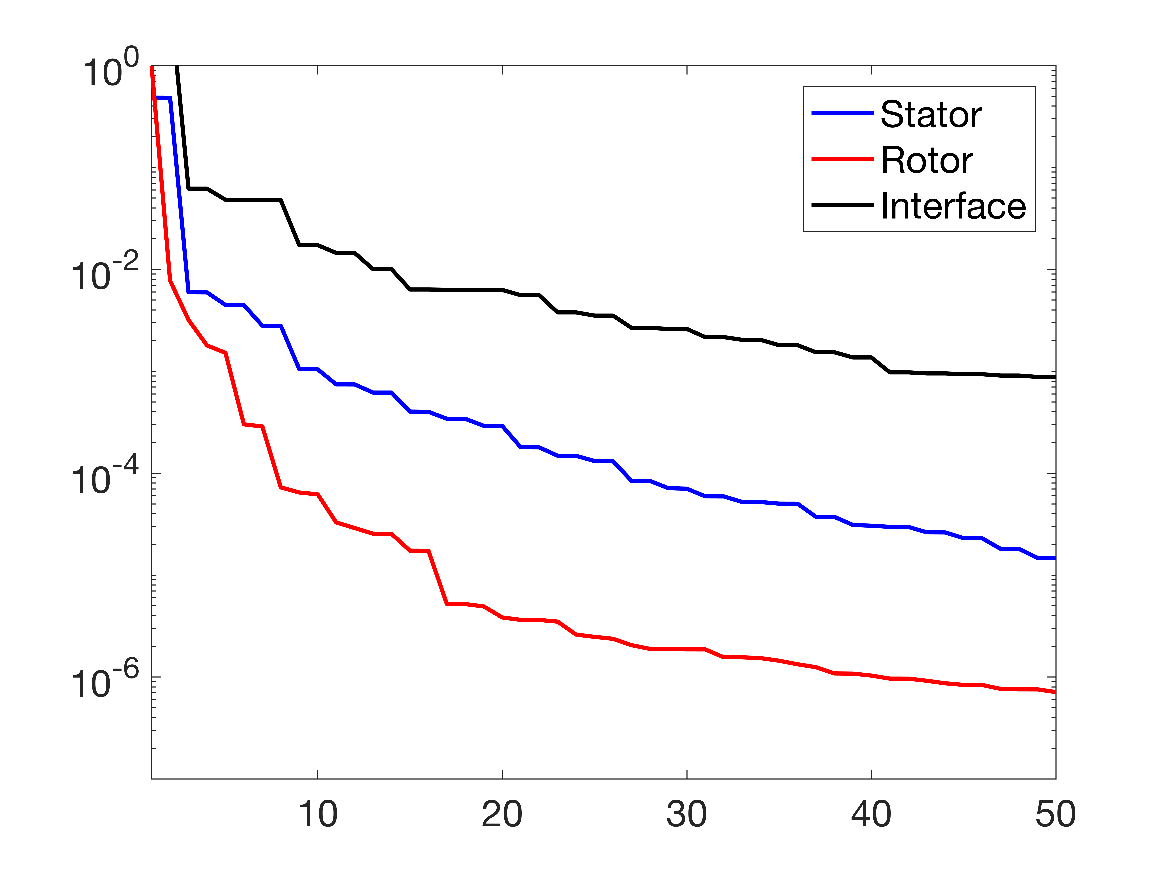}
        \\
          (a) \emph{sym}
        & (b) \emph{rot}
        \\[6pt]
          \includegraphics[width=58mm]{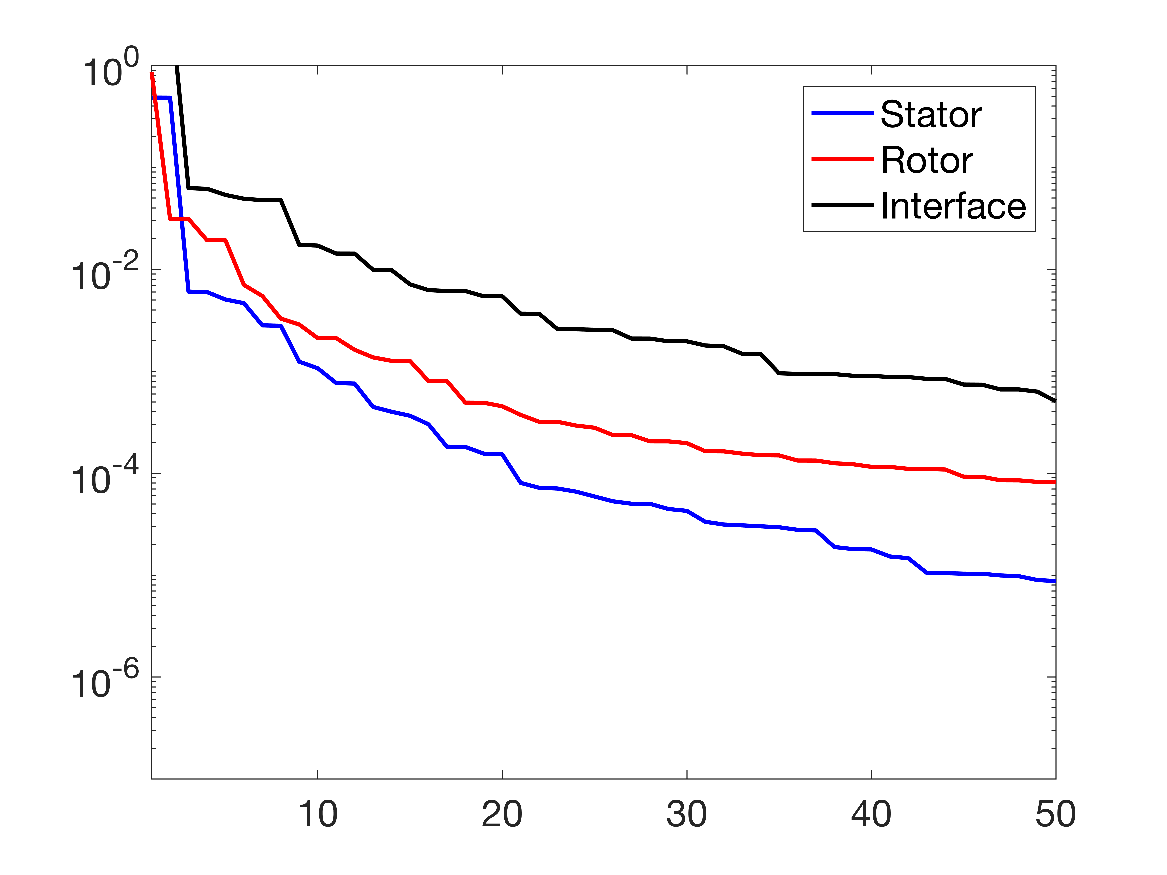}
        & \includegraphics[width=58mm]{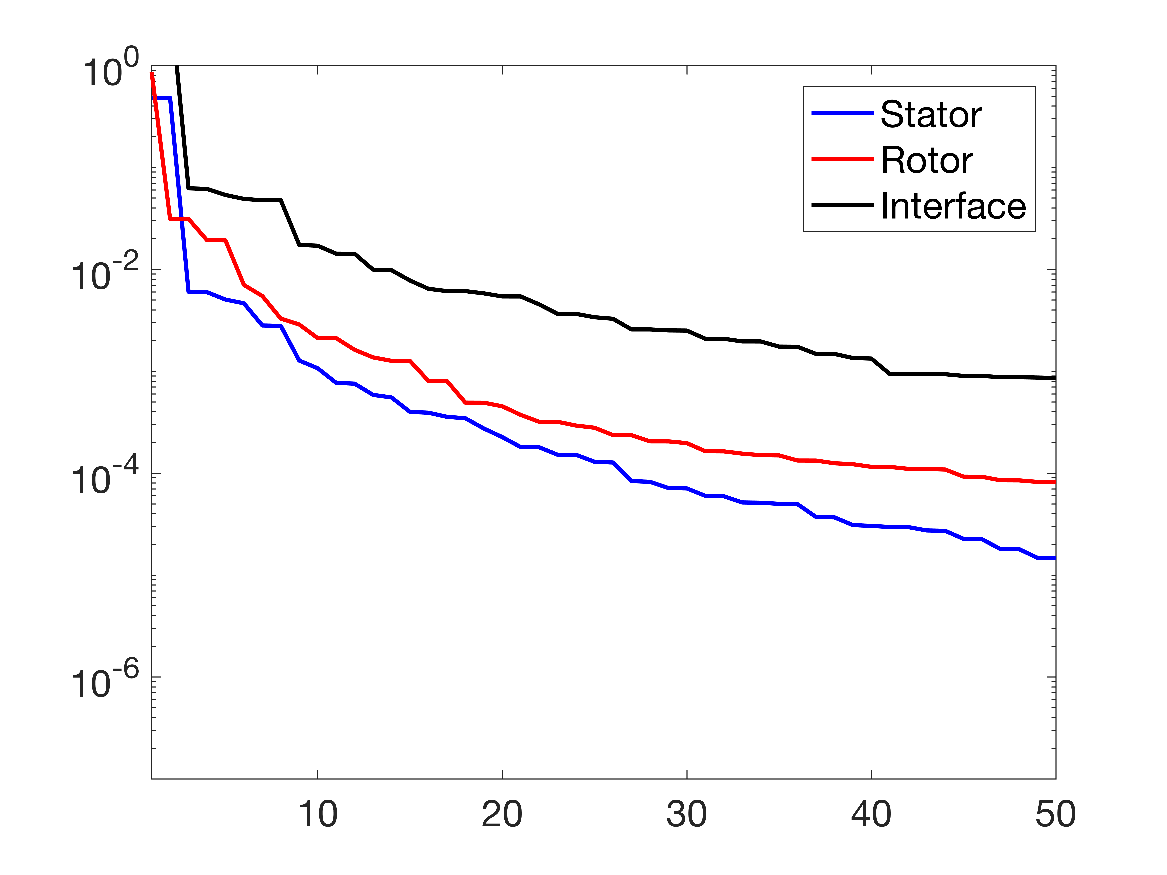}
        \\
          (c) \emph{stat}
        & (d) \emph{rot\_stat}
    \end{tabular}
    \caption{Decay of the normalized eigenvalues for the different settings.}
    \label{fig:eig}
\end{figure}

To start we will have a look at the performance of the POD method. For this we
compute a full revolution for each of the four setting and have a look at the
decay of the eigenvalues. Let us recall that a fast decay is essential for the
POD method to perform well. In Figure~\ref{fig:eig} the decay of the normalized
eigenvalues is shown. As can be seen the eigenvalues decay very fast for rotor
and stator. Only for the interface the decay is much slower. This underlines
the decision to not perform a model order reduction for the interface.
Additionally it can be observed that a perturbation in the magnetic field angle
$\phi$ has less influence on the decay of the eigenvalues than the perturbation
in the length $\ell$ of one stator tooth. While the perturbation in $\phi$
causes a slight change in the decay of the eigenvalues related to the stator,
the perturbation in $\ell$ dramatically influences the eigenvalues related to
the rotor.

\begin{figure}
\centering
 \includegraphics[width=35mm]{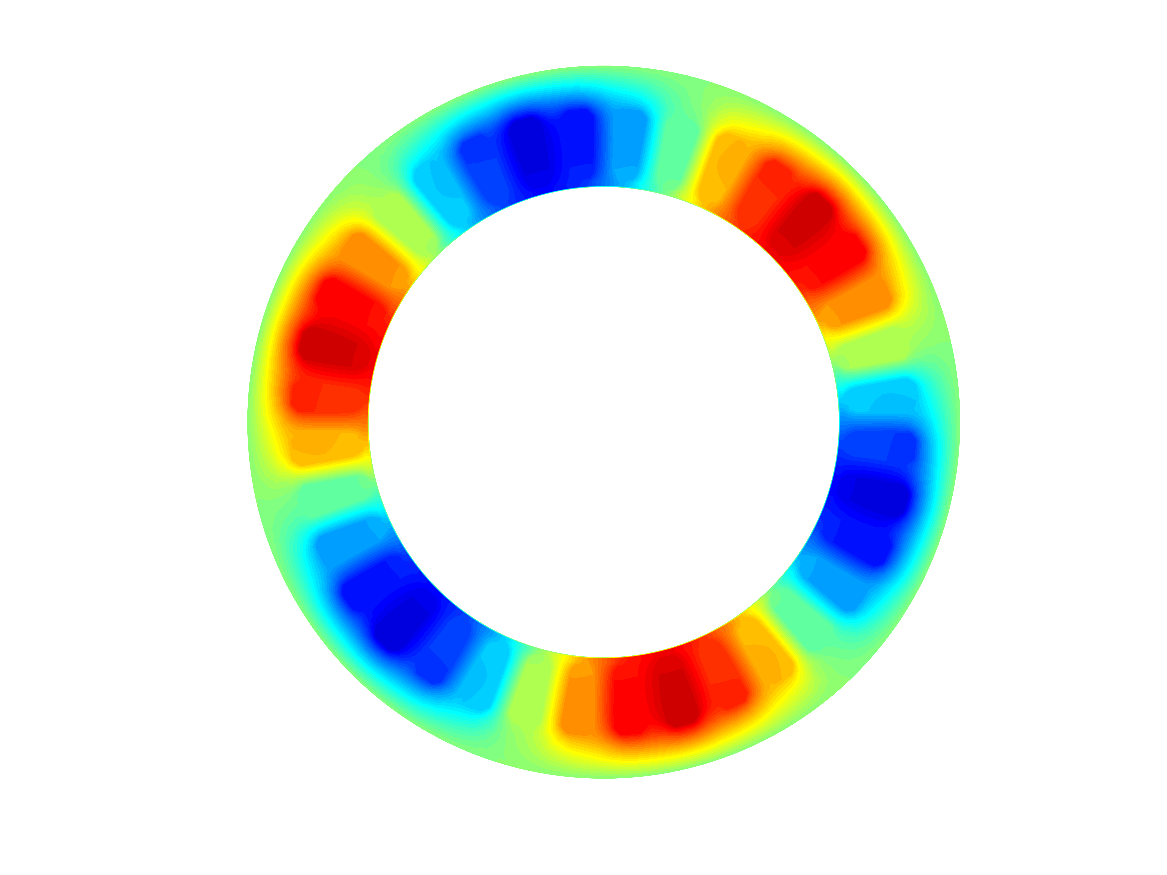}
 \includegraphics[width=35mm]{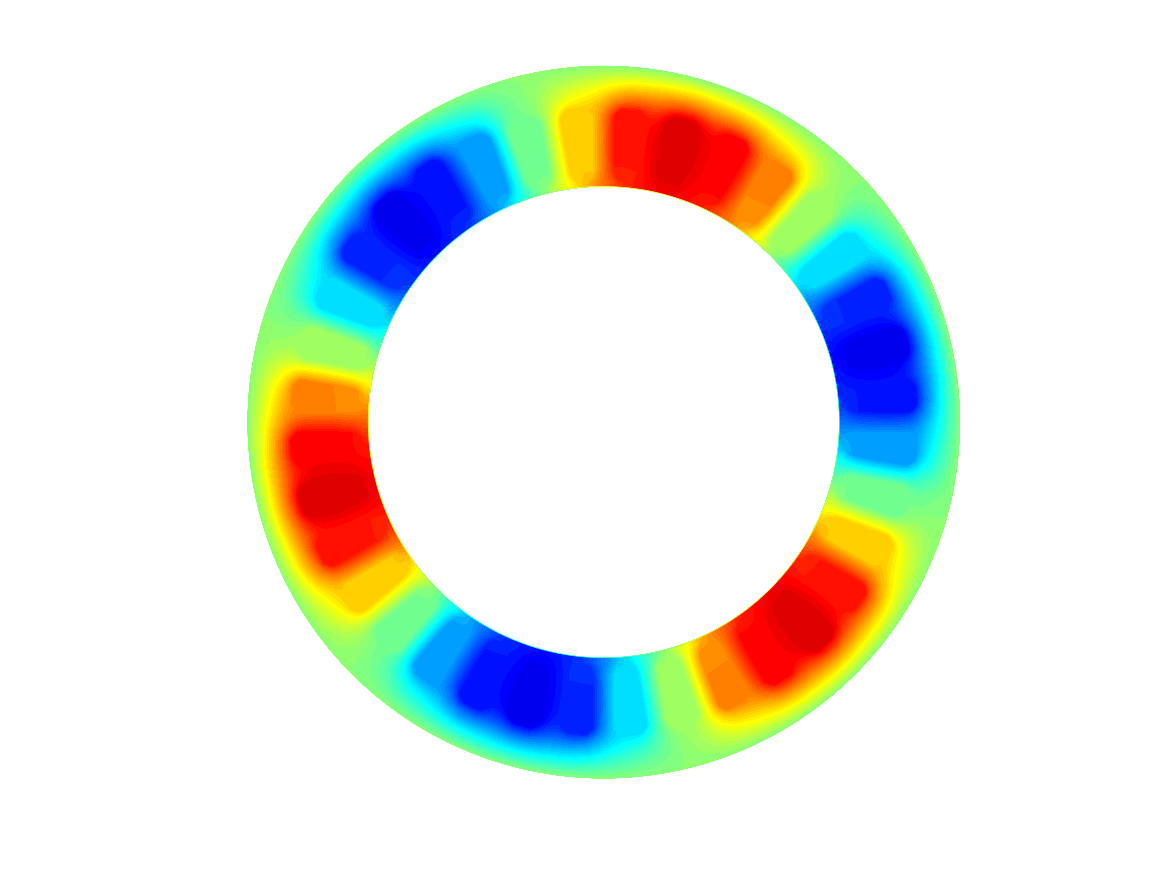}
 \hspace{-4mm}\includegraphics[width=45mm]{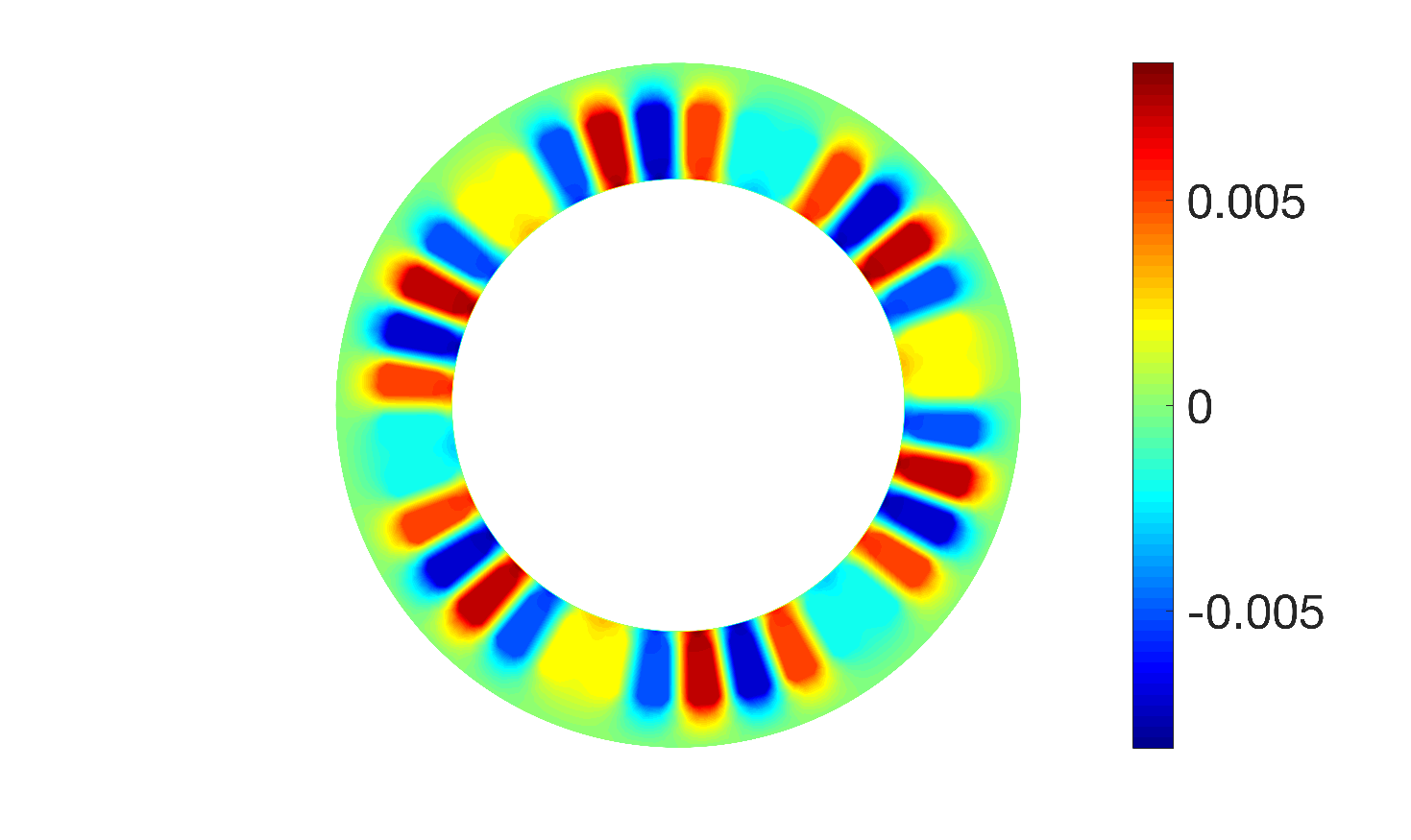}
 \includegraphics[width=35mm]{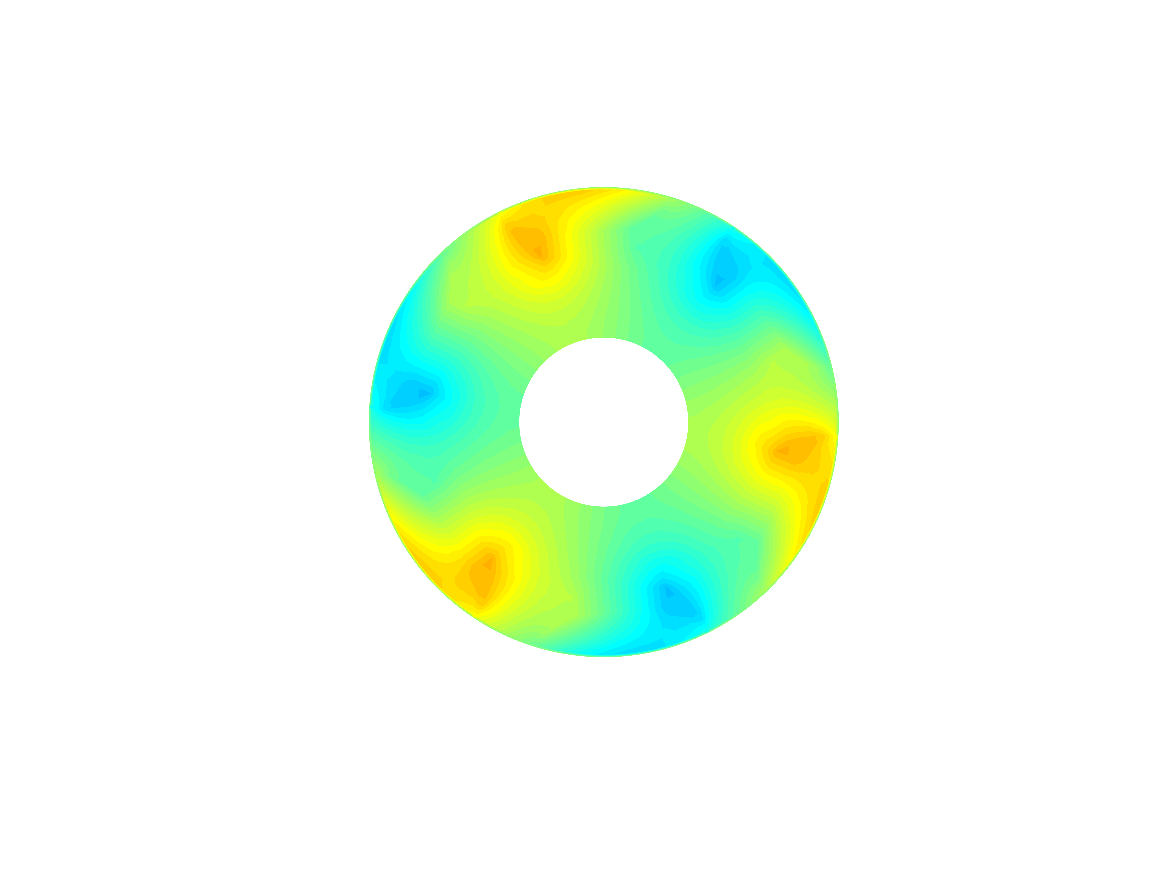}
 \includegraphics[width=35mm]{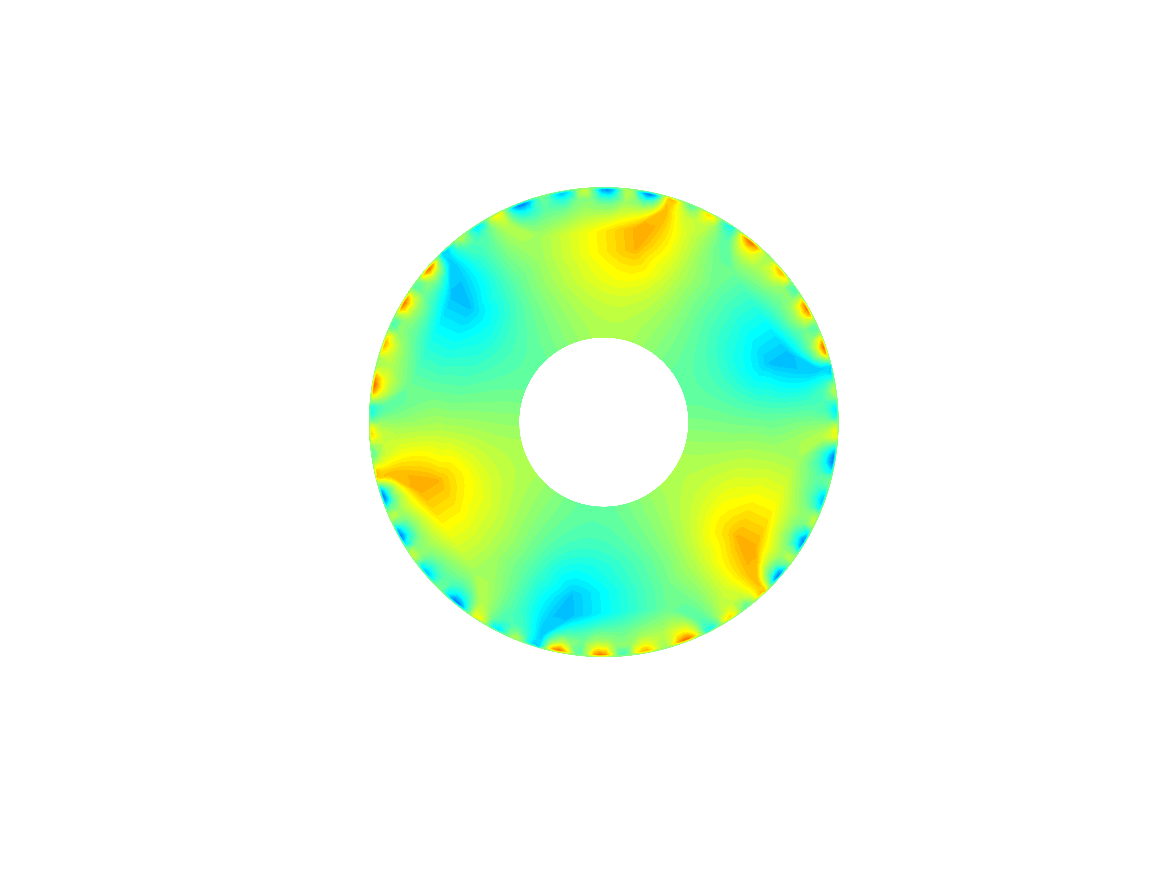}
 \hspace{-4mm}\includegraphics[width=45mm]{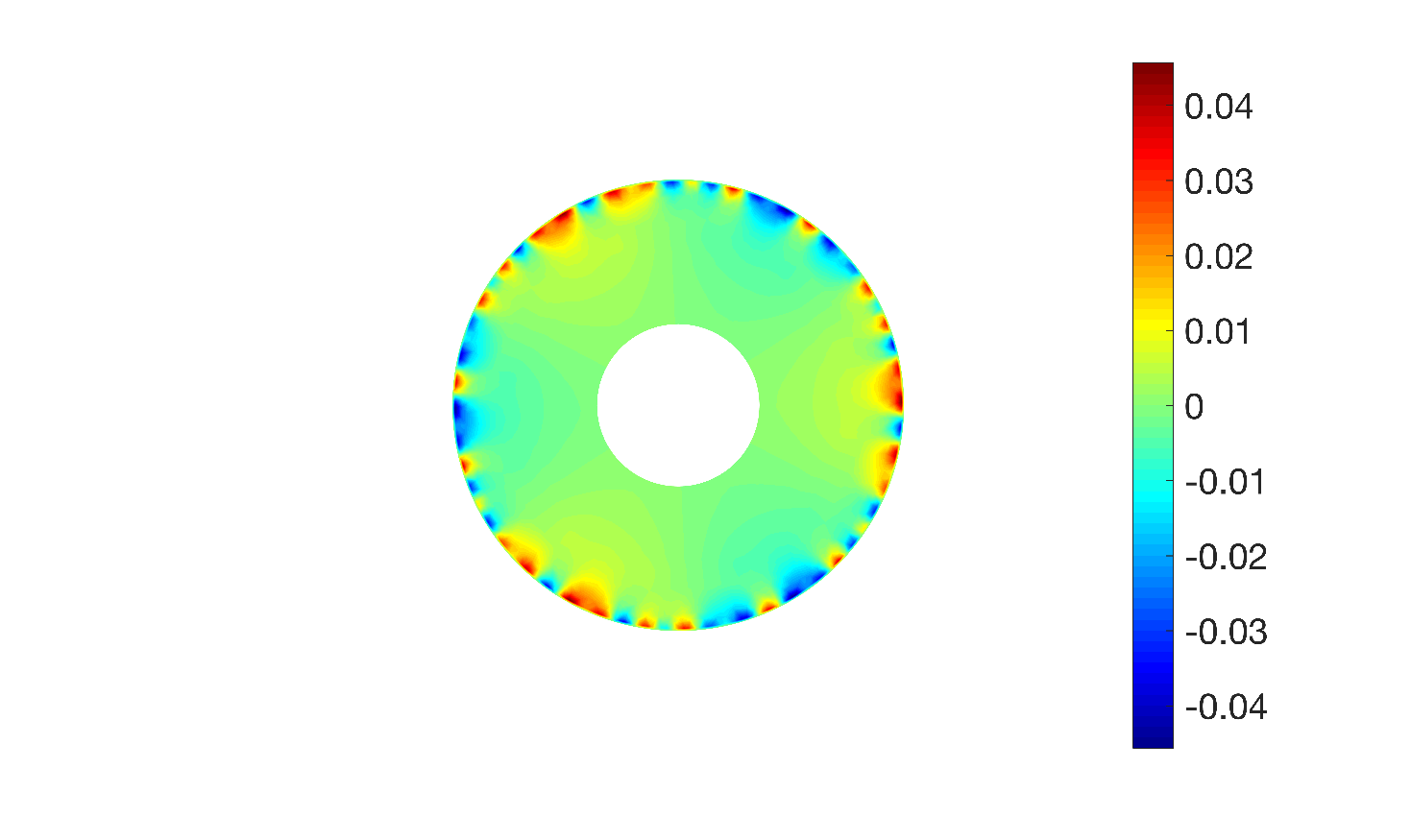}
 \caption{First three POD basis vectors for the stator (top) and rotor (bottom) for setting \emph{sym}.}\label{fig:pod_basis}
\end{figure}

\begin{figure}
\centering
 \includegraphics[width=35mm]{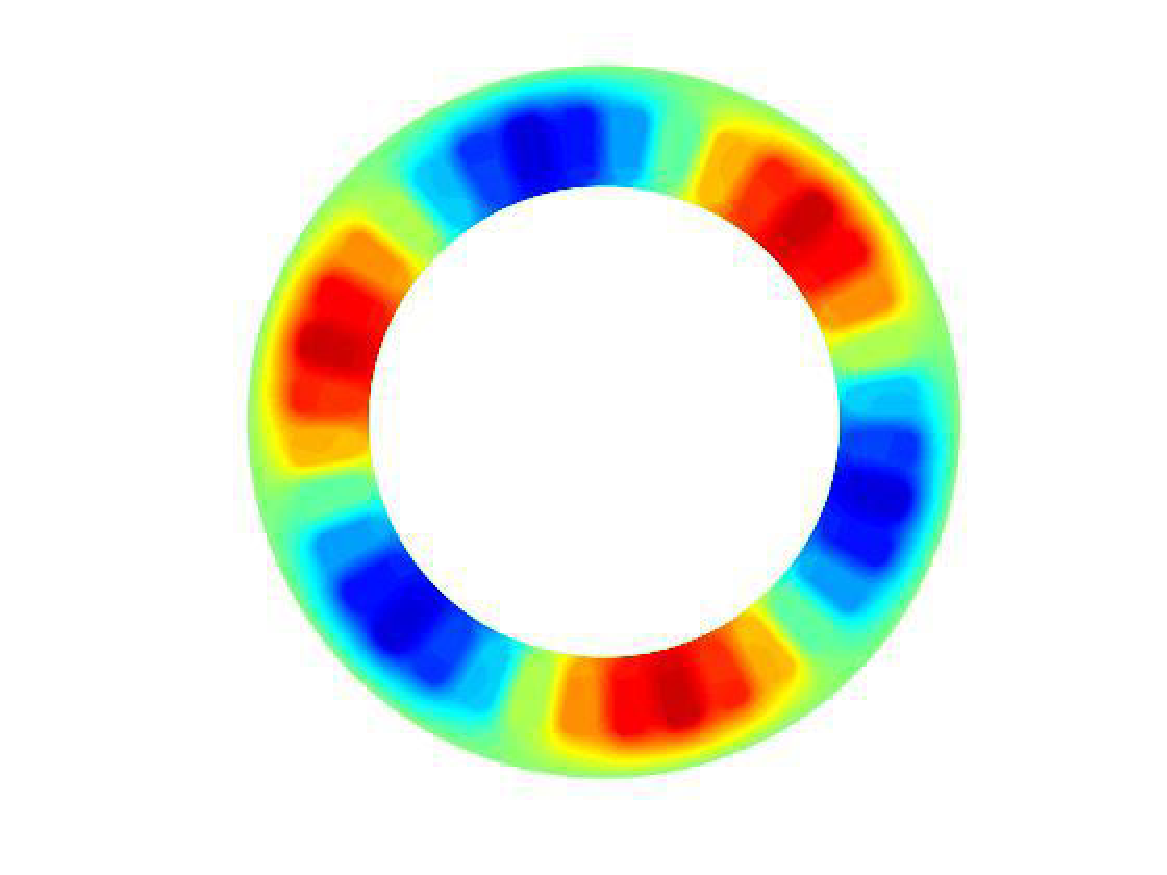}
 \includegraphics[width=35mm]{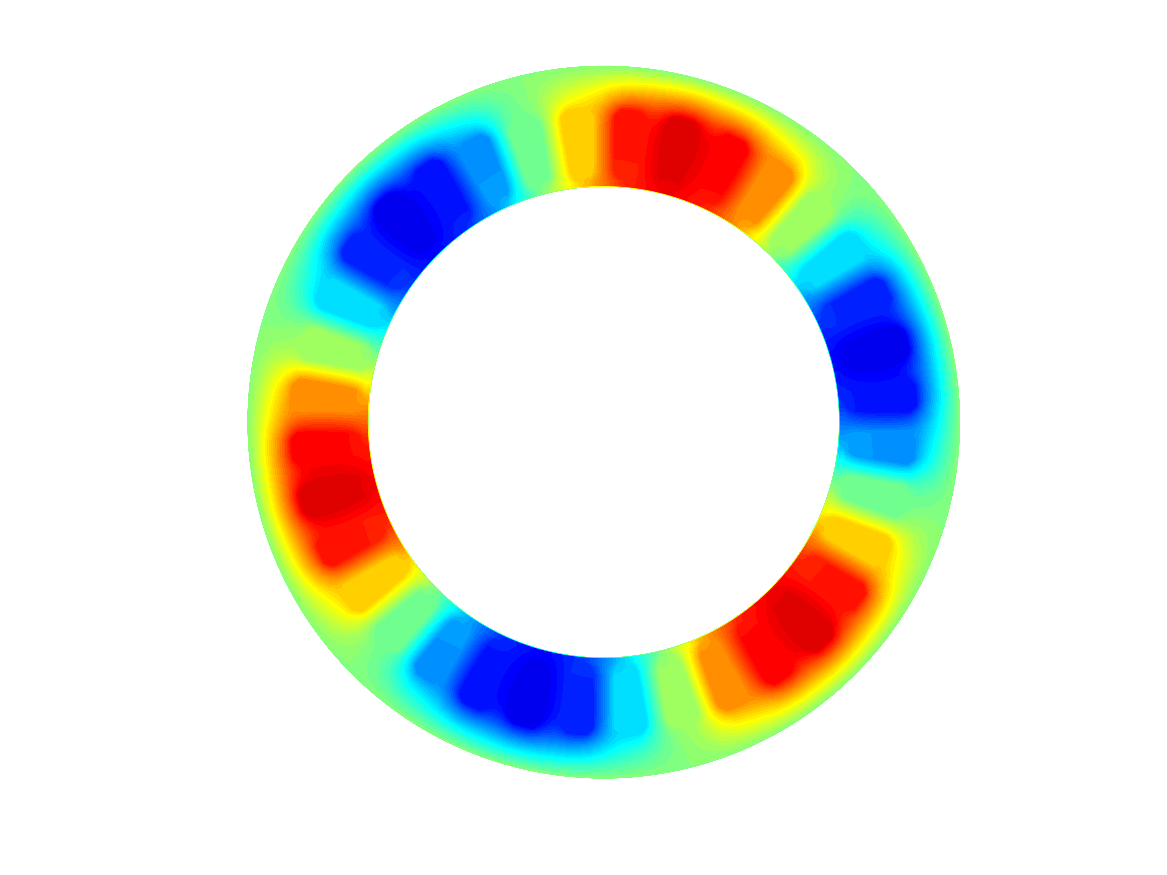}
 \hspace{-4mm}\includegraphics[width=45mm]{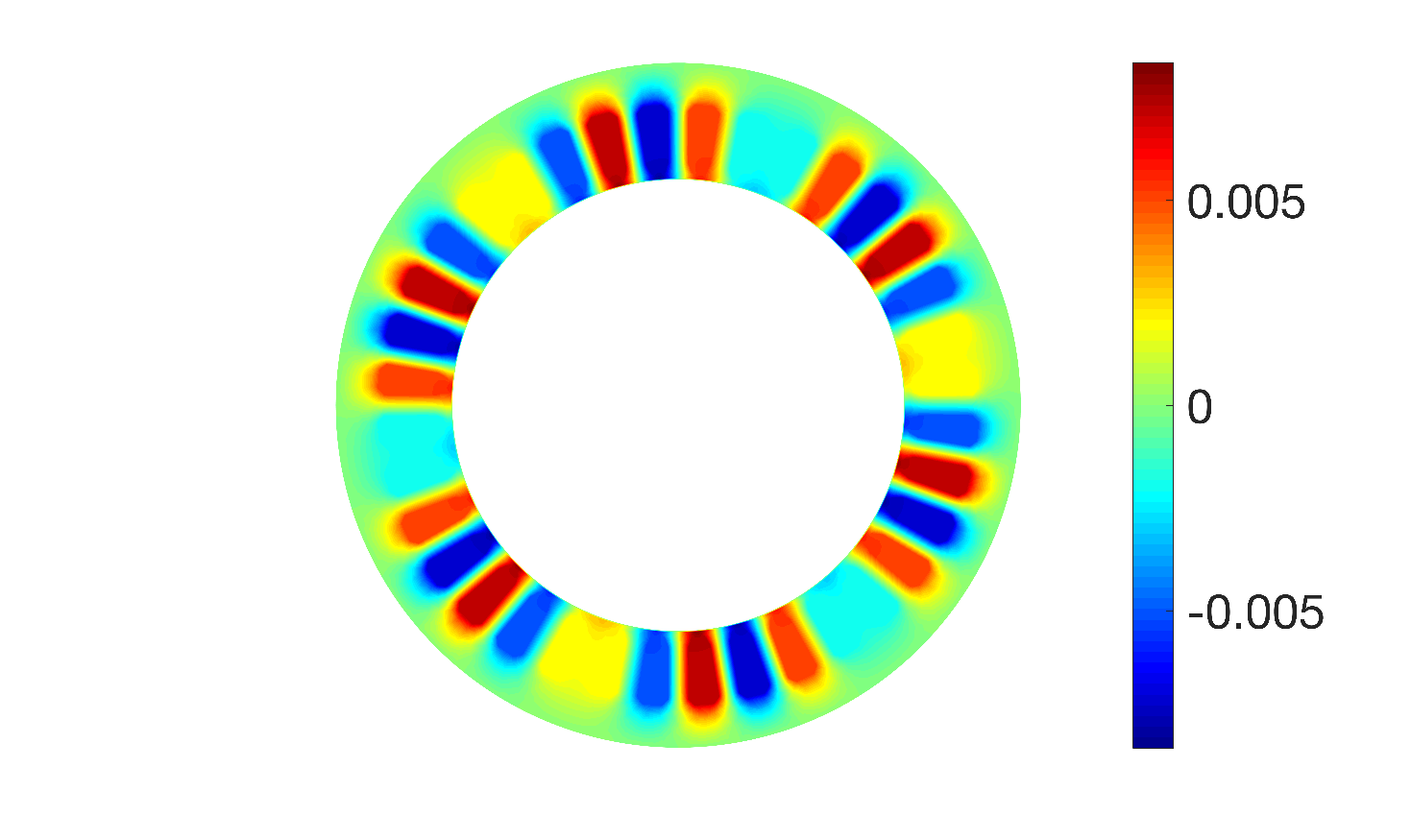}
 \includegraphics[width=35mm]{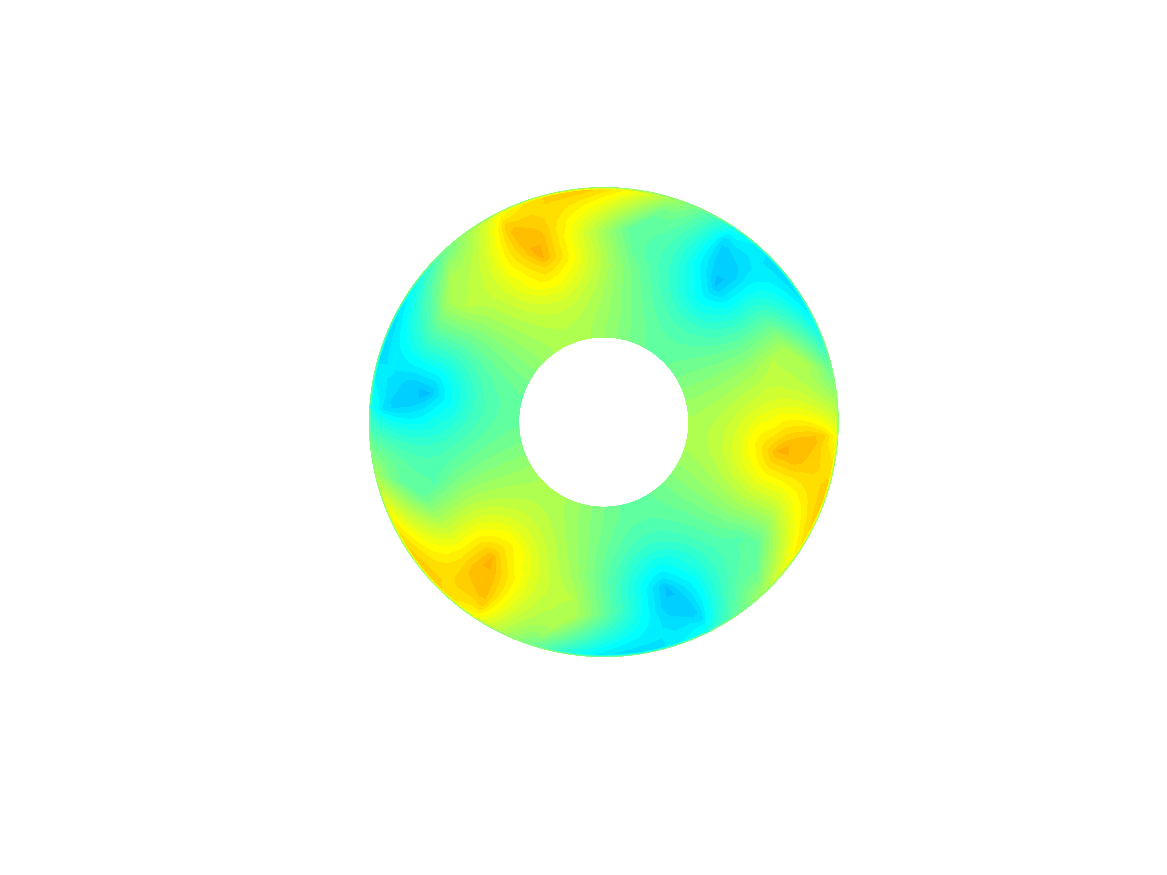}
 \includegraphics[width=35mm]{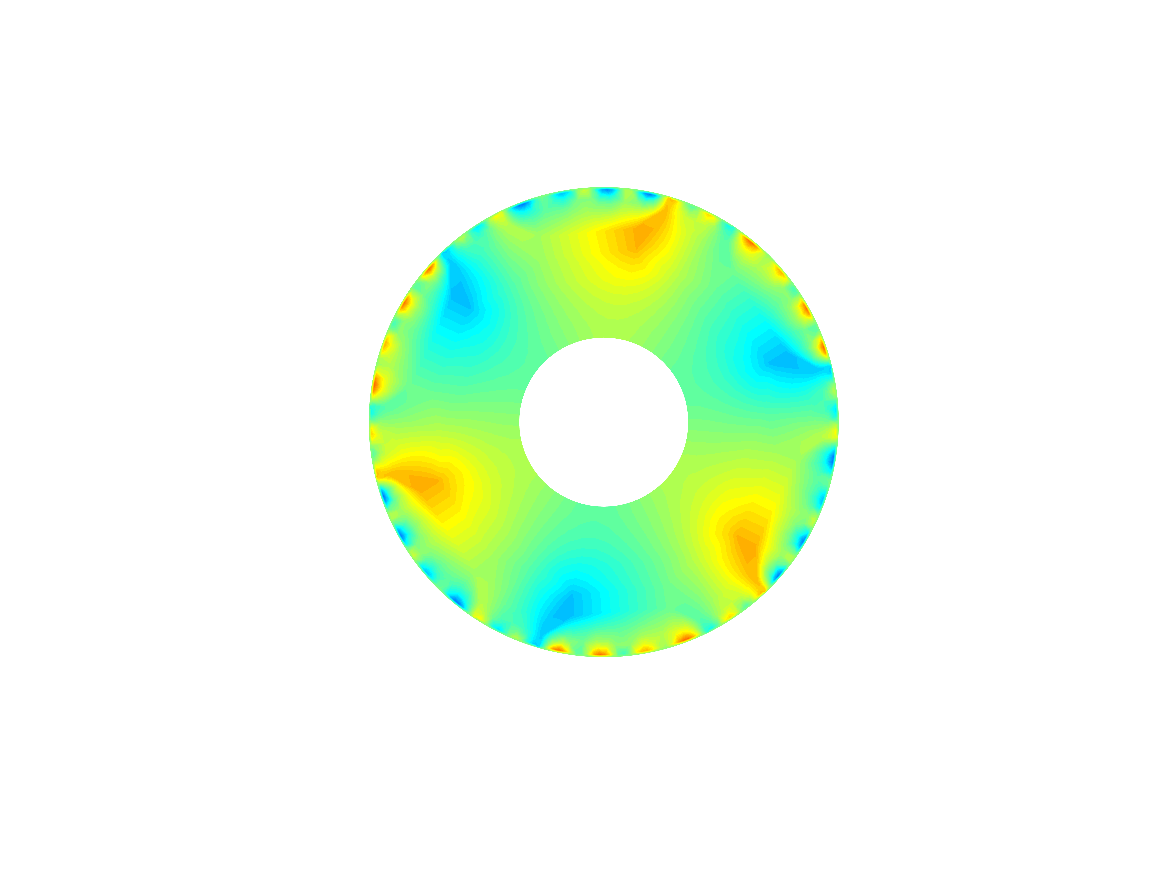}
 \hspace{-4mm}\includegraphics[width=45mm]{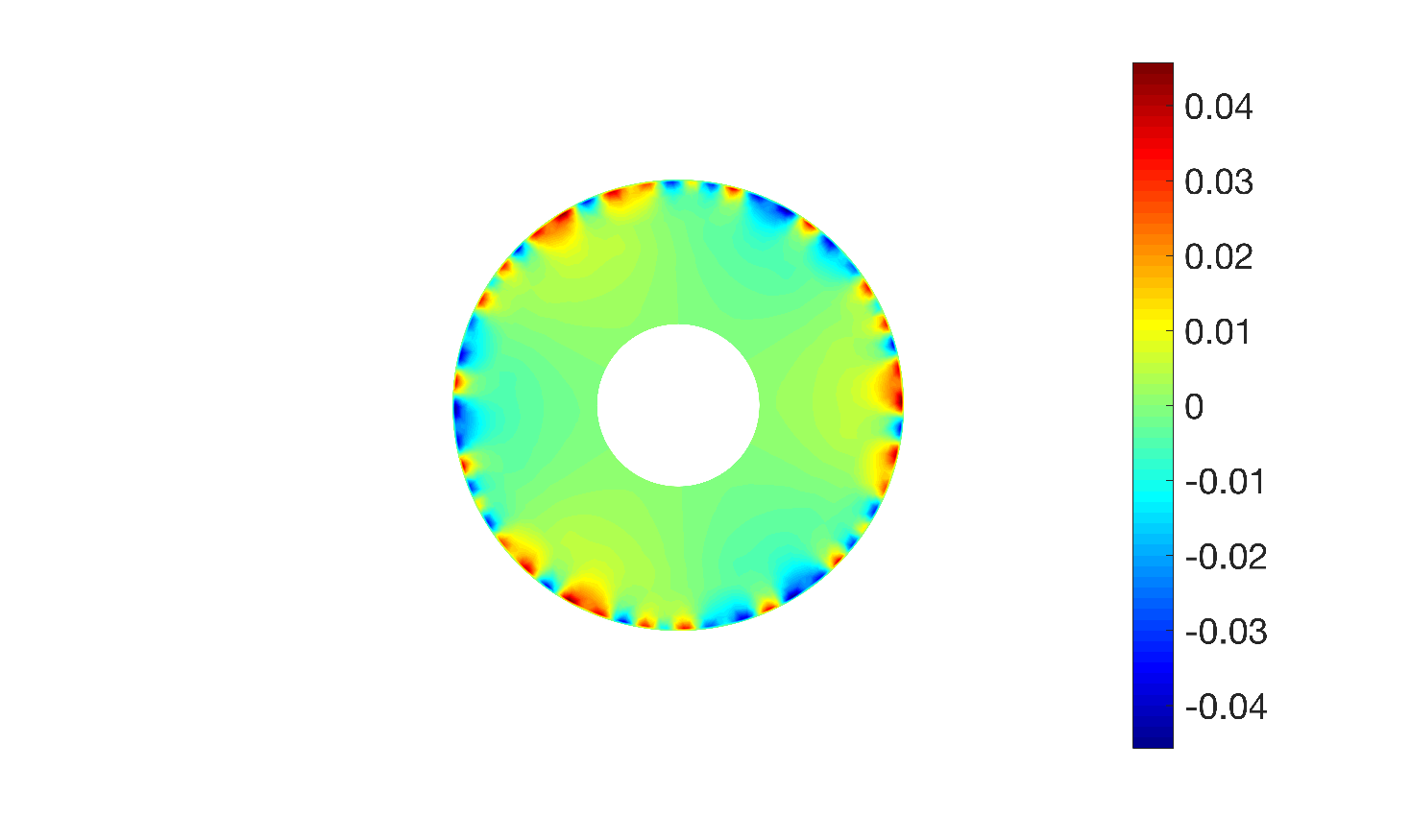}
 \caption{First three POD basis vectors for the stator (top) and rotor (bottom) for setting \emph{rot}.}\label{fig:pod_basis_r}
\end{figure}

\begin{figure}
\centering
 \includegraphics[width=35mm]{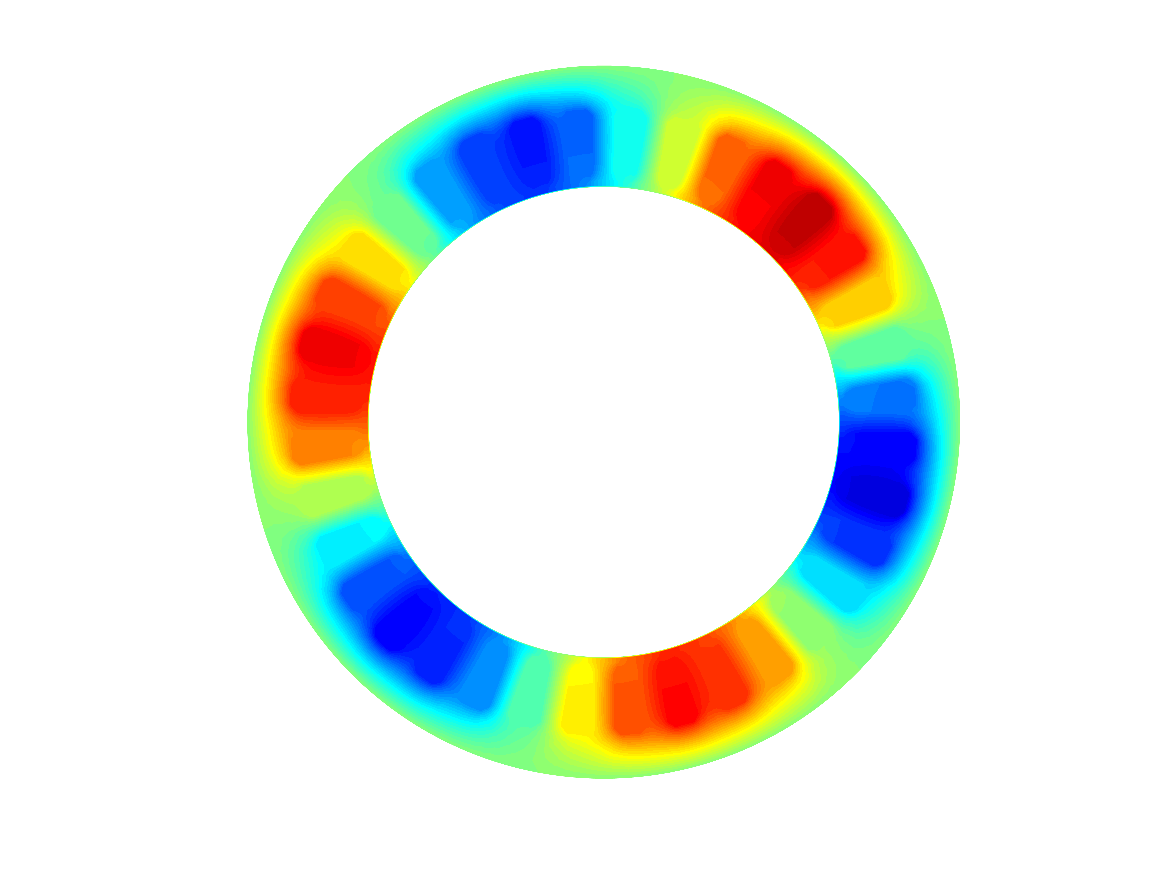}
 \includegraphics[width=35mm]{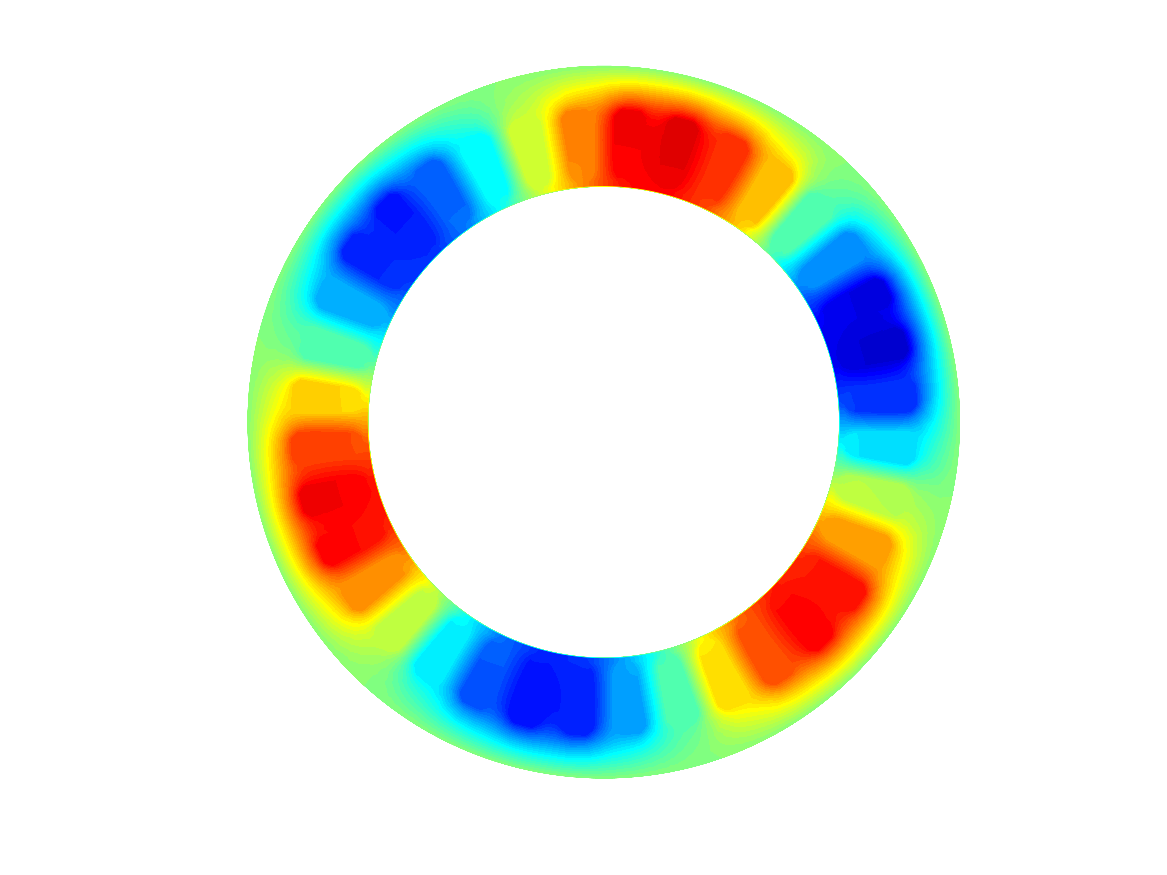}
 \hspace{-4mm}\includegraphics[width=45mm]{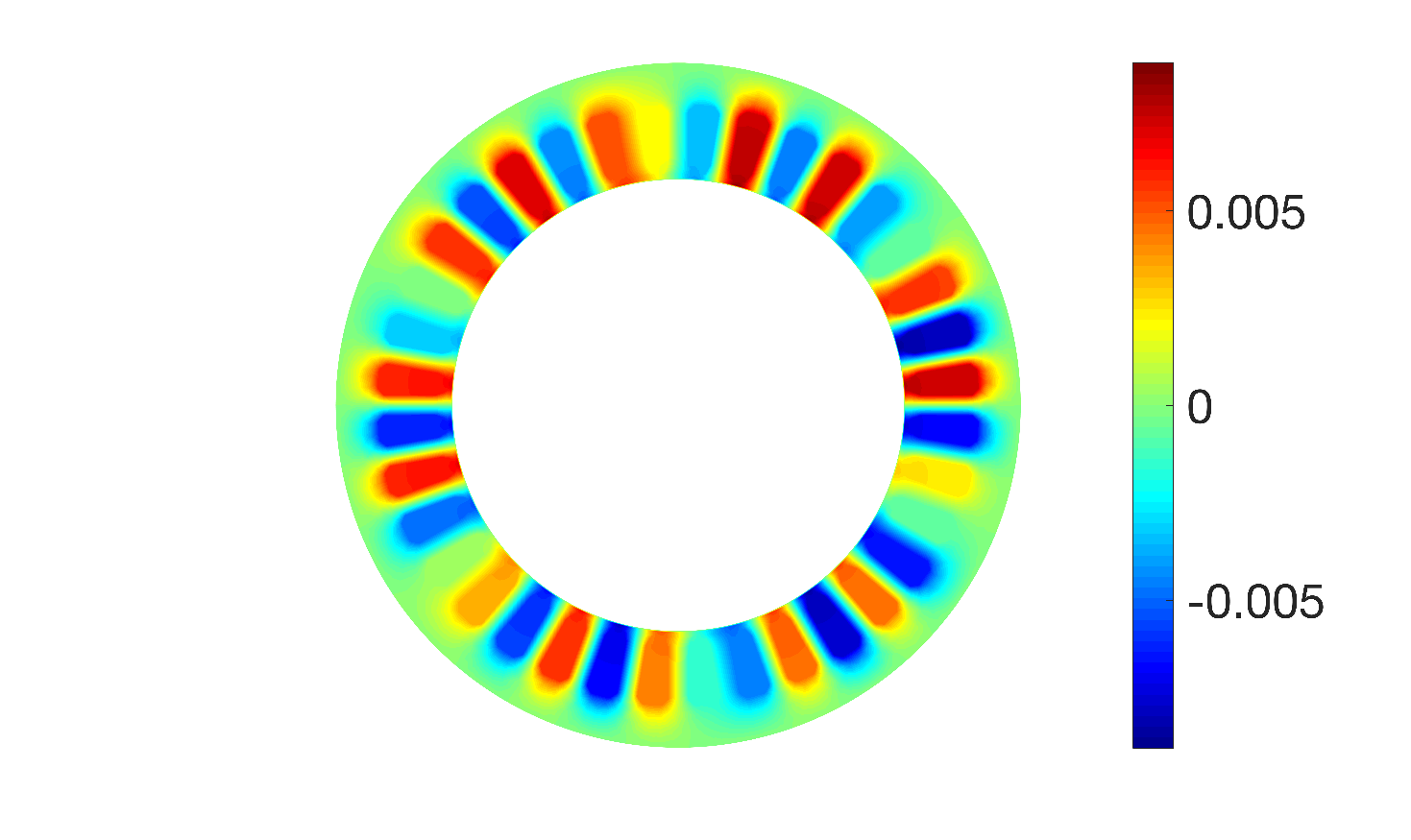}
 \includegraphics[width=35mm]{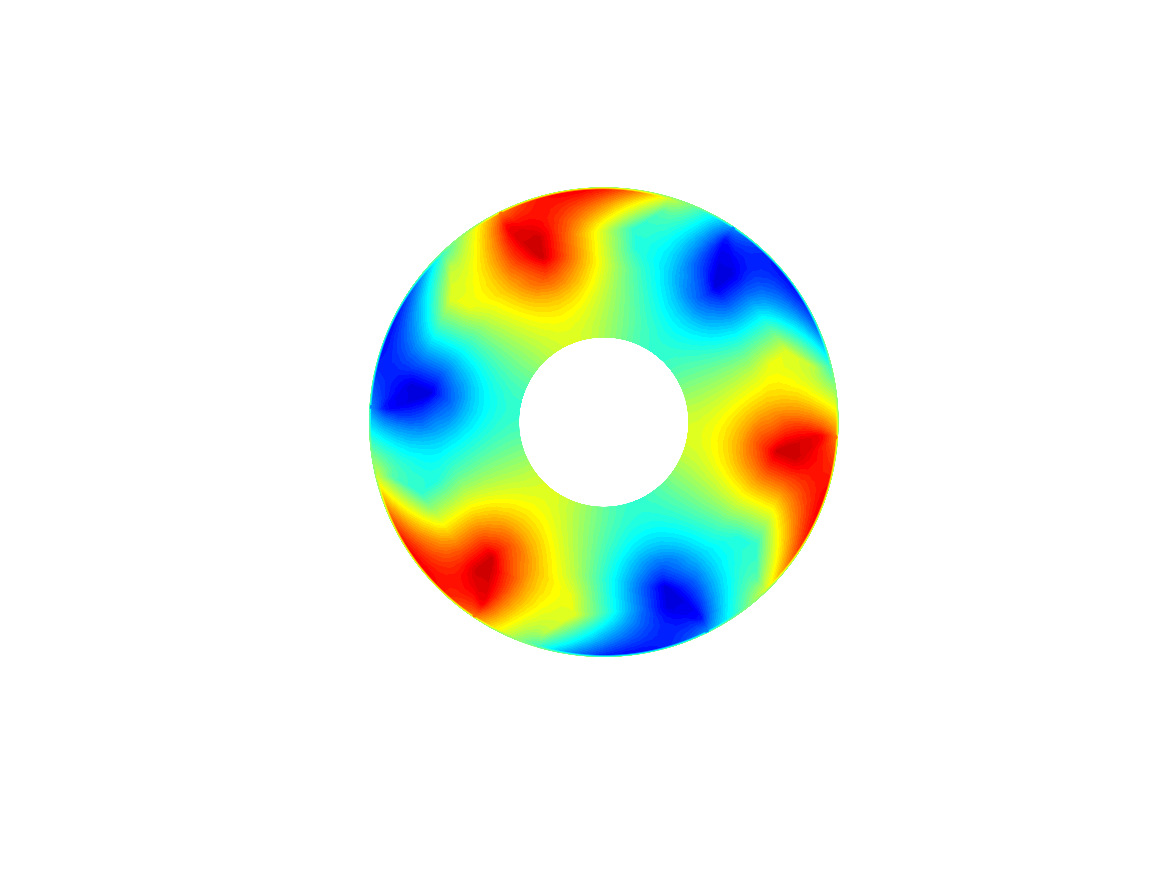}
 \includegraphics[width=35mm]{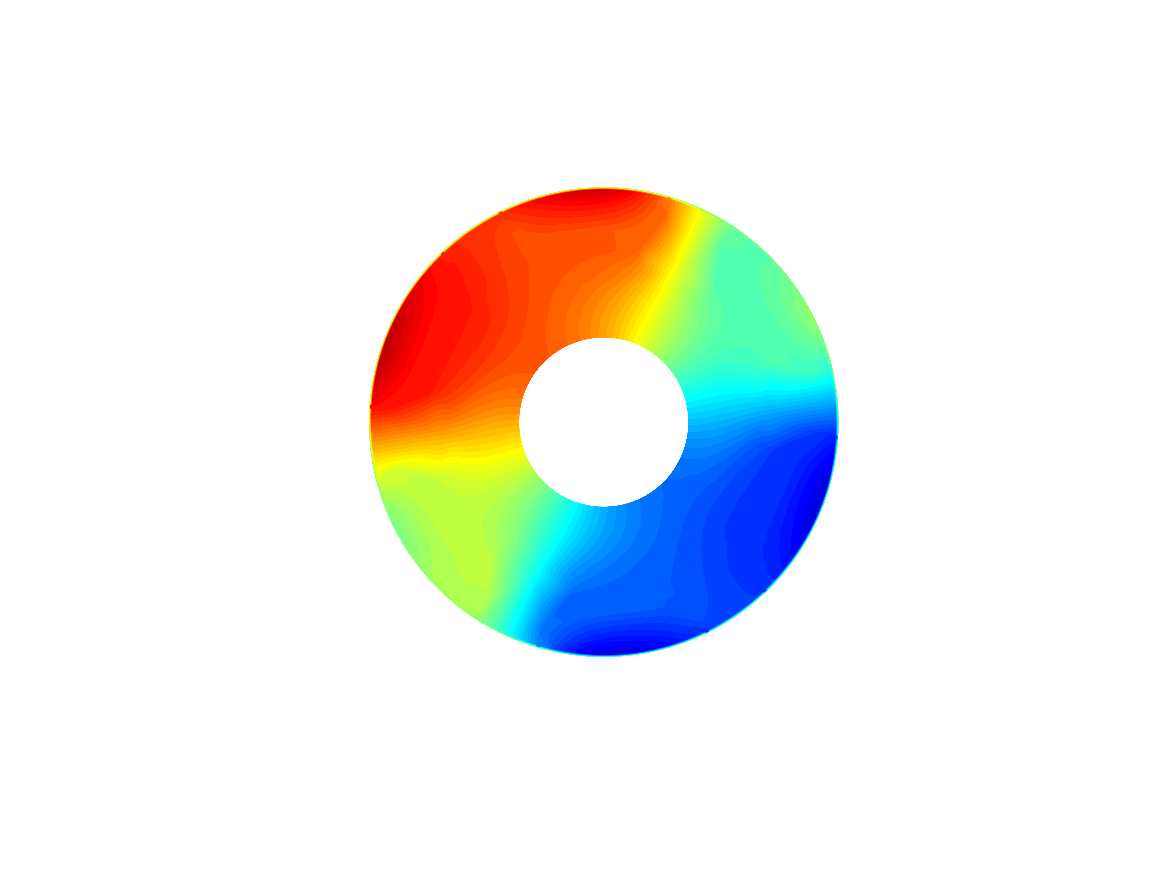}
 \hspace{-4mm}\includegraphics[width=45mm]{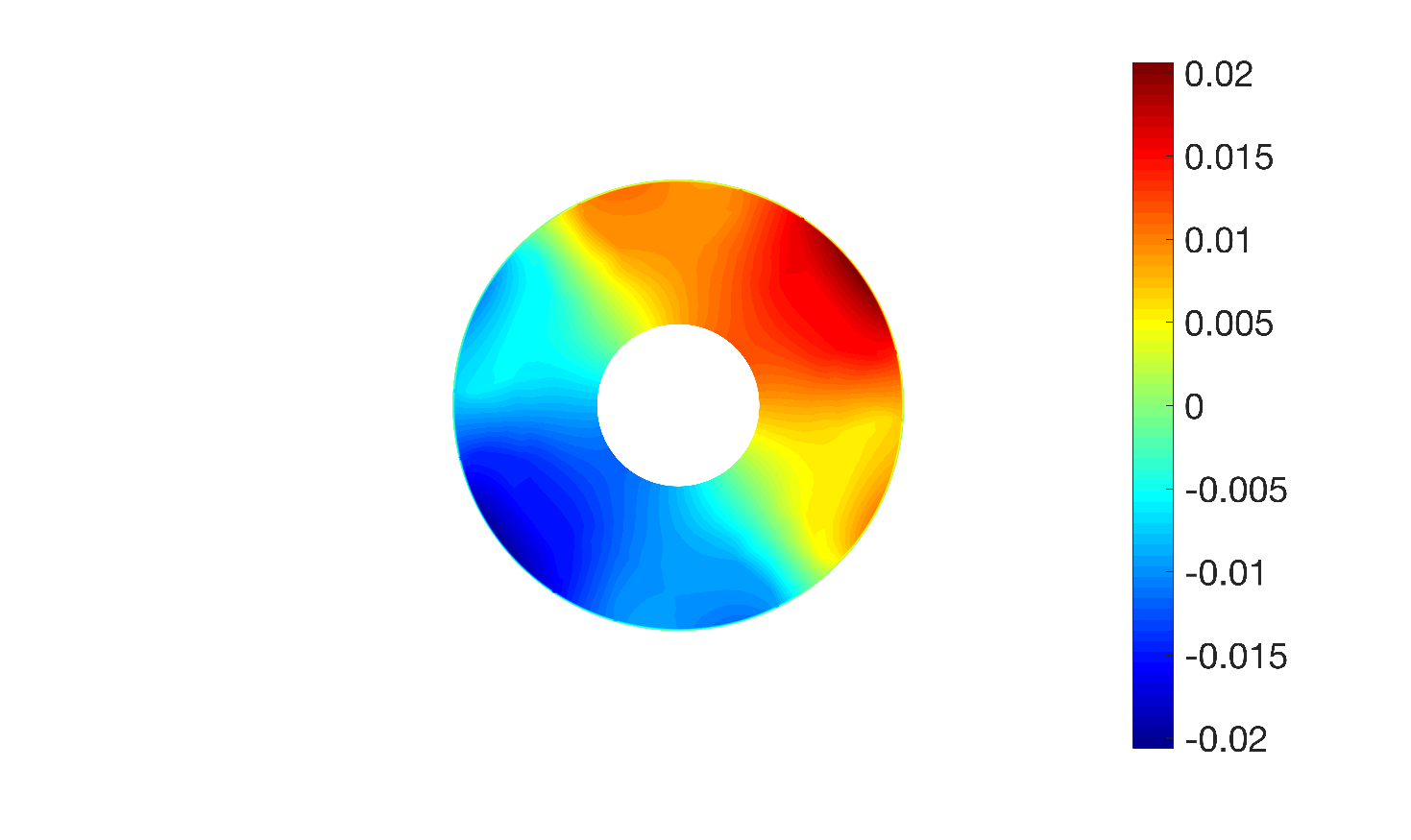}
 \caption{First three POD basis vectors for the stator (top) and rotor (bottom) for setting \emph{stat}.}\label{fig:pod_basis_s}
\end{figure}

\begin{figure}
\centering
 \includegraphics[width=35mm]{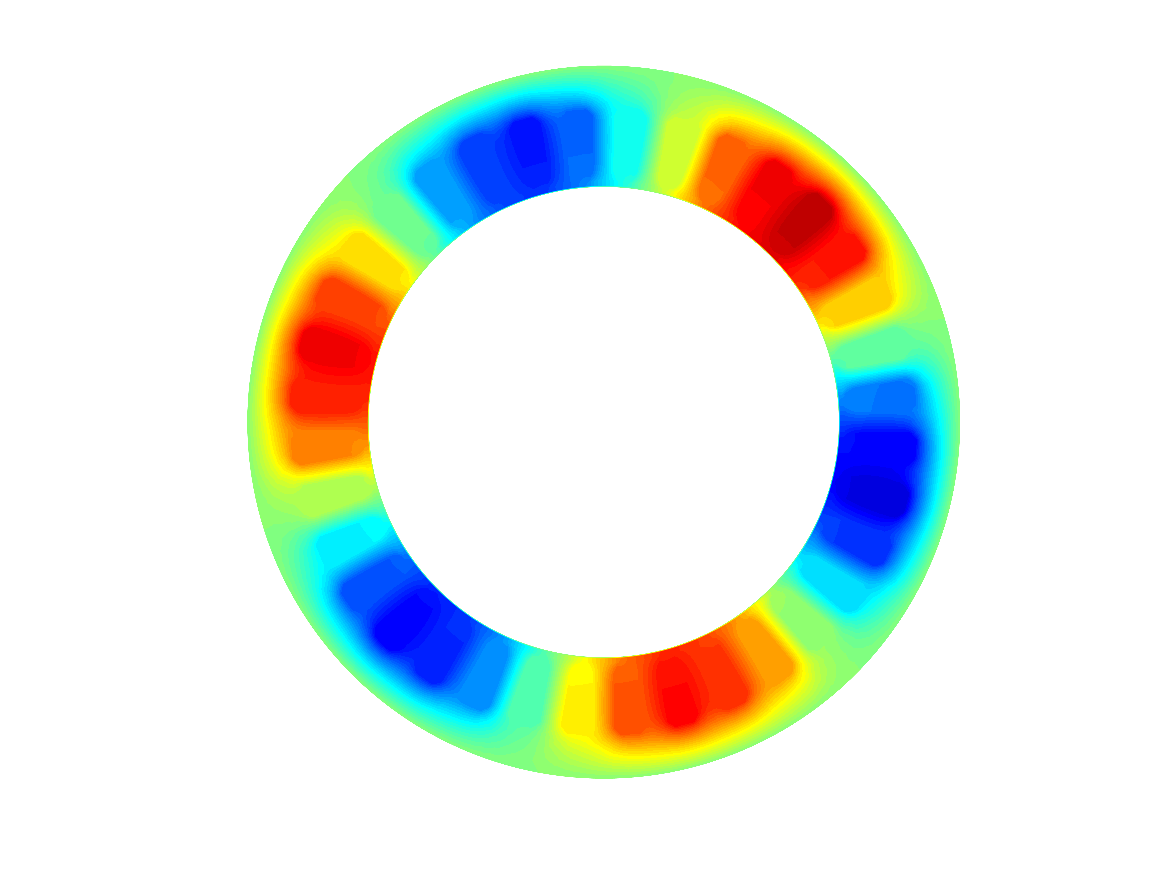}
 \includegraphics[width=35mm]{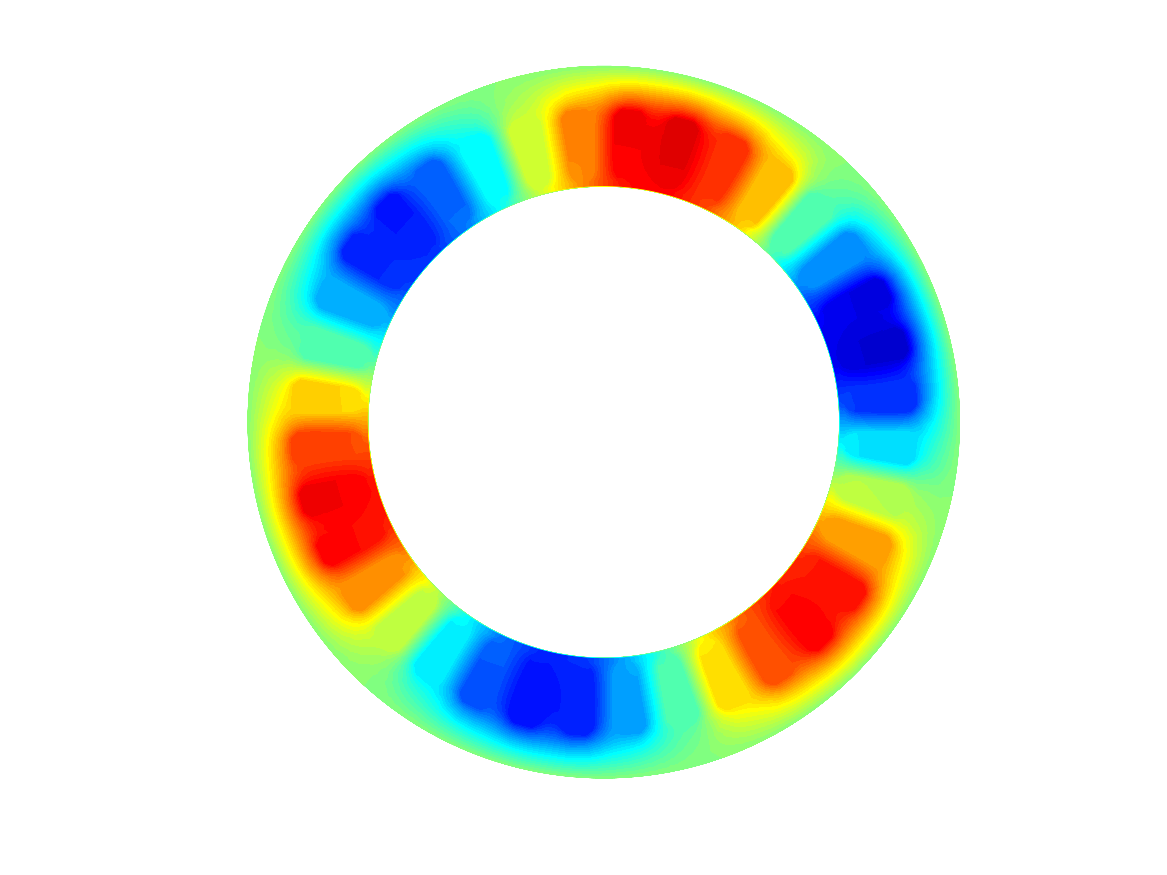}
 \hspace{-4mm}\includegraphics[width=45mm]{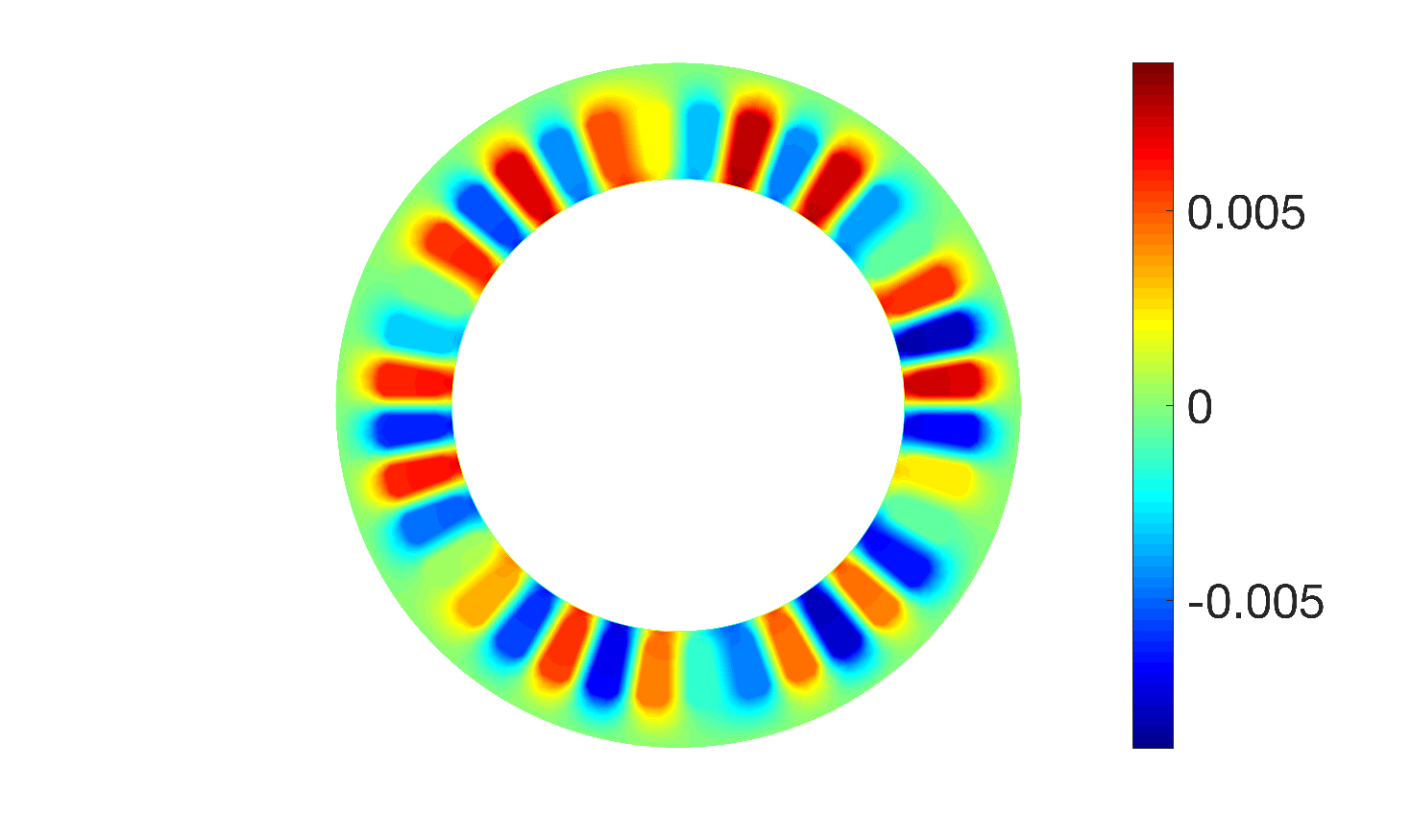}
 \includegraphics[width=35mm]{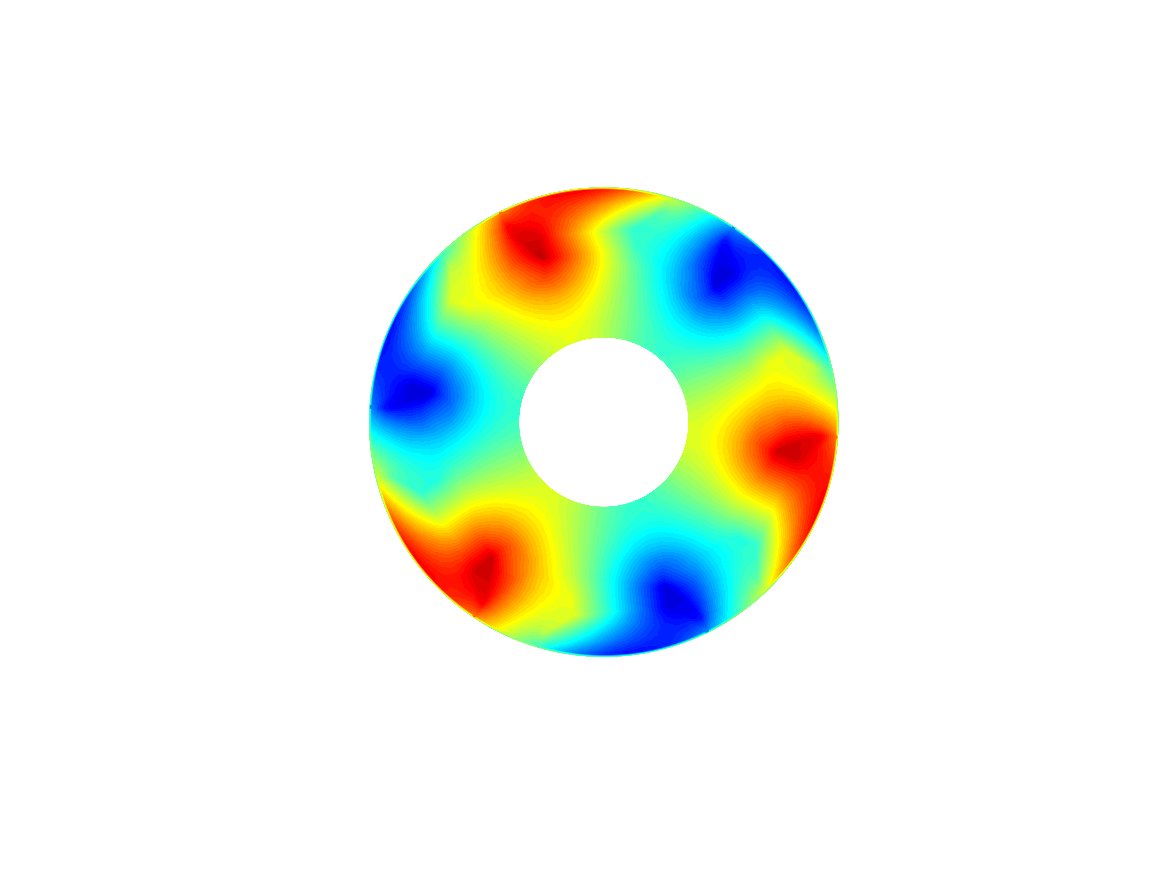}
 \includegraphics[width=35mm]{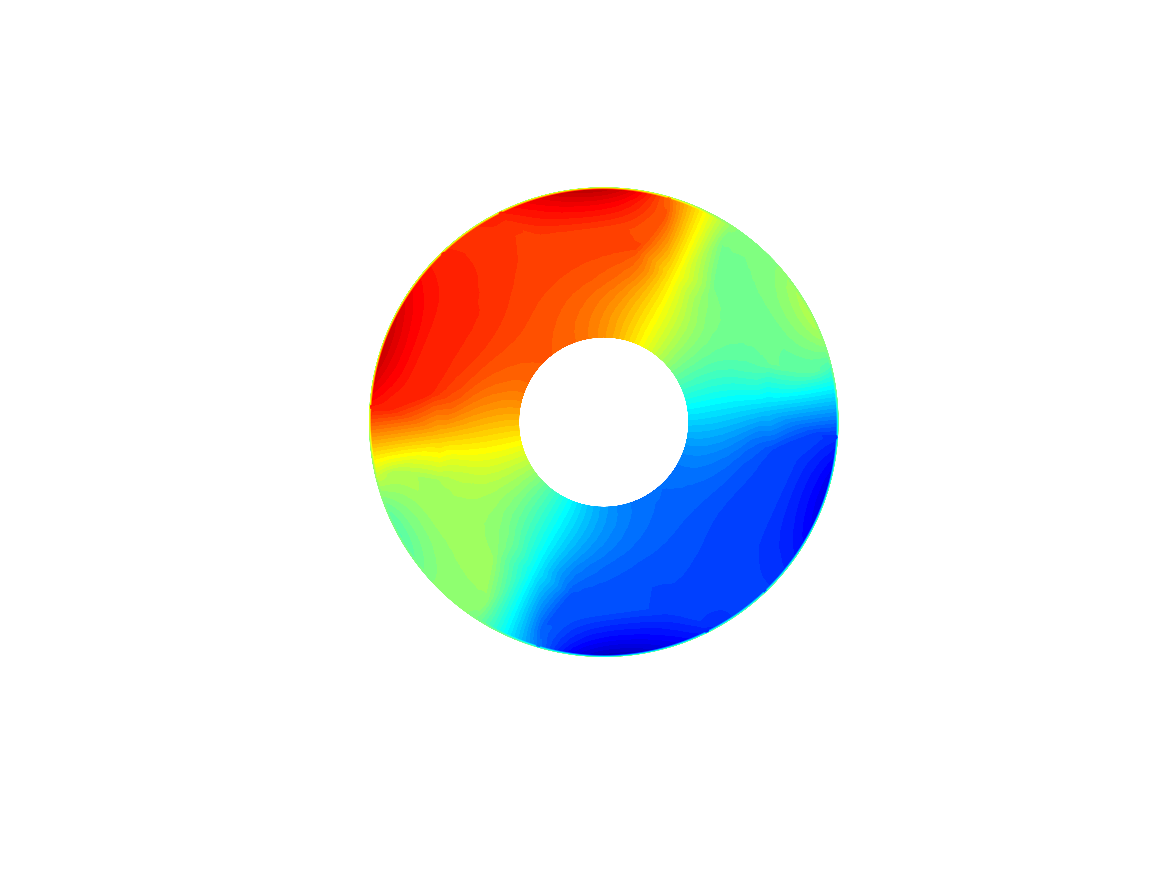}
 \hspace{-4mm}\includegraphics[width=45mm]{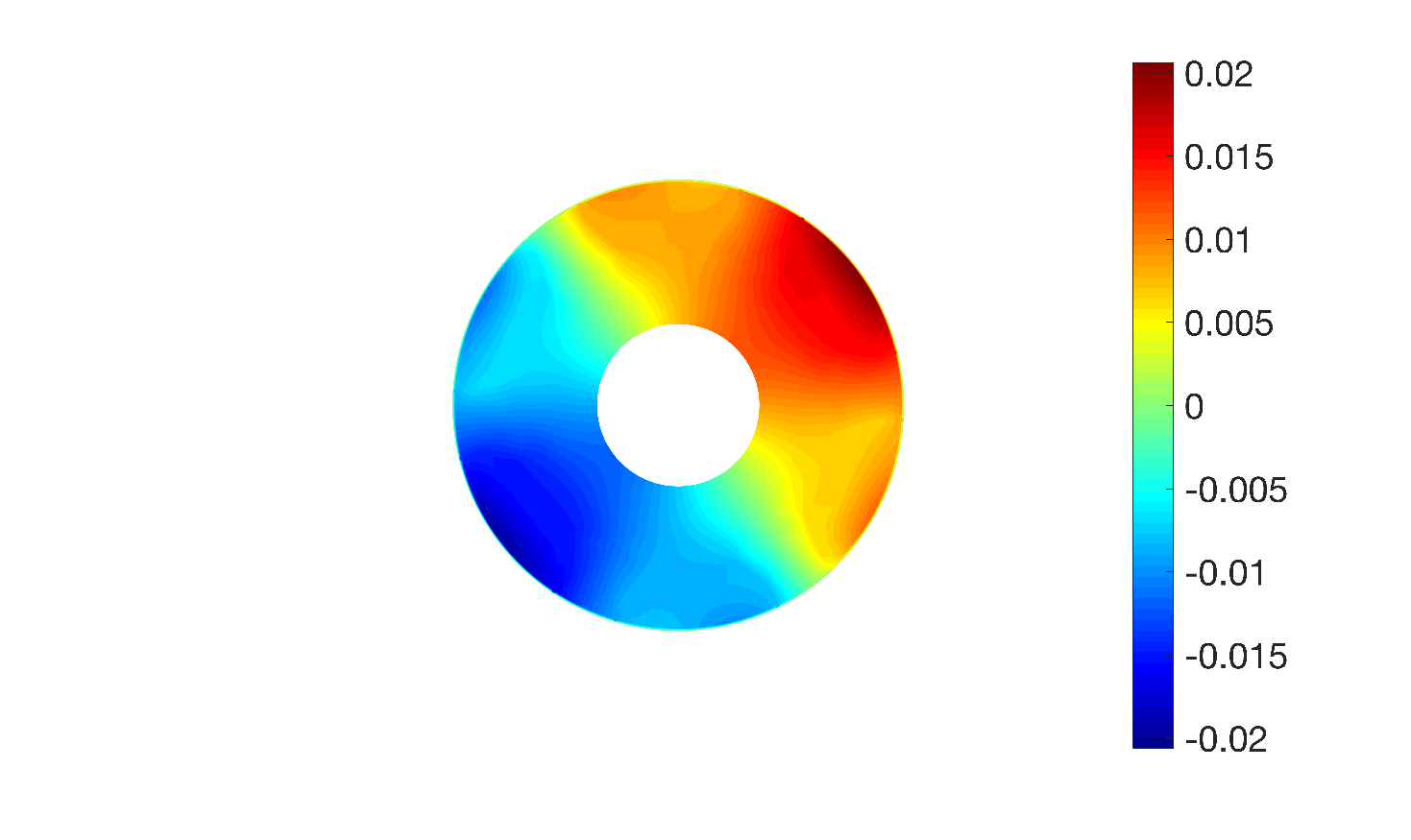}
 \caption{First three POD basis vectors for the stator (top) and rotor (bottom) for setting \emph{rot\_stat}.}\label{fig:pod_basis_sr}
\end{figure}

Next we depict the first three POD basis vectors for each setting
(Figure~\ref{fig:pod_basis}-\ref{fig:pod_basis_sr}). For the setting \emph{sym}
and \emph{rot} very similar basis vectors are obtained. The last basis vector
for both settings is zero on wide areas and only adds contributions at the
interface (green is zero). In the settings \emph{stat} and \emph{rot\_stat}
this is very different. There the third basis vector for the rotor still
contributes a lot of information to the reduced order model. What was already
observed in the decay of the eigenvalues was verified again in the POD basis
vectors. The settings \emph{sym} and \emph{rot} exhibit similar behaviors as
well as \emph{stat} and \emph{rot\_stat}.

From the decay of the eigenvalues it can be expected that we will be able to
generate reduced order models of very low dimension. The plots of POD basis
verify this for the \emph{sym} and \emph{rot} since already the third basis
vector is almost zero.

Next we have a look at the performance of the presented adaptive algorithm. For
this let us introduce the index sets $\mathcal K_i$. We will look at two
settings. To be able to differentiate the two settings we will call sets
$\mathcal M$ and $\mathcal K$. In the first we choose the sets by dividing the
full rotation into six sections (6 poles). within each section we have a
hierarchy of sets for refinement. Hence the sets can be written as follows
\begin{equation*} \mathcal K_{ij} = 150(i-1) + [j, j+12, j+24, \ldots ],\quad i
= 1,\ldots,6\quad\mbox{and}\quad j=1,\ldots,12. \end{equation*} Note that the
values for $\mathcal K_{ij}$ are always less than $150i$. Using these sets we
always compute snapshots only related to one pole. If the error is still to
large in the particular pole a refinement is achieved.

\begin{figure}
    \begin{tabular}{cc}
           \includegraphics[width=58mm]{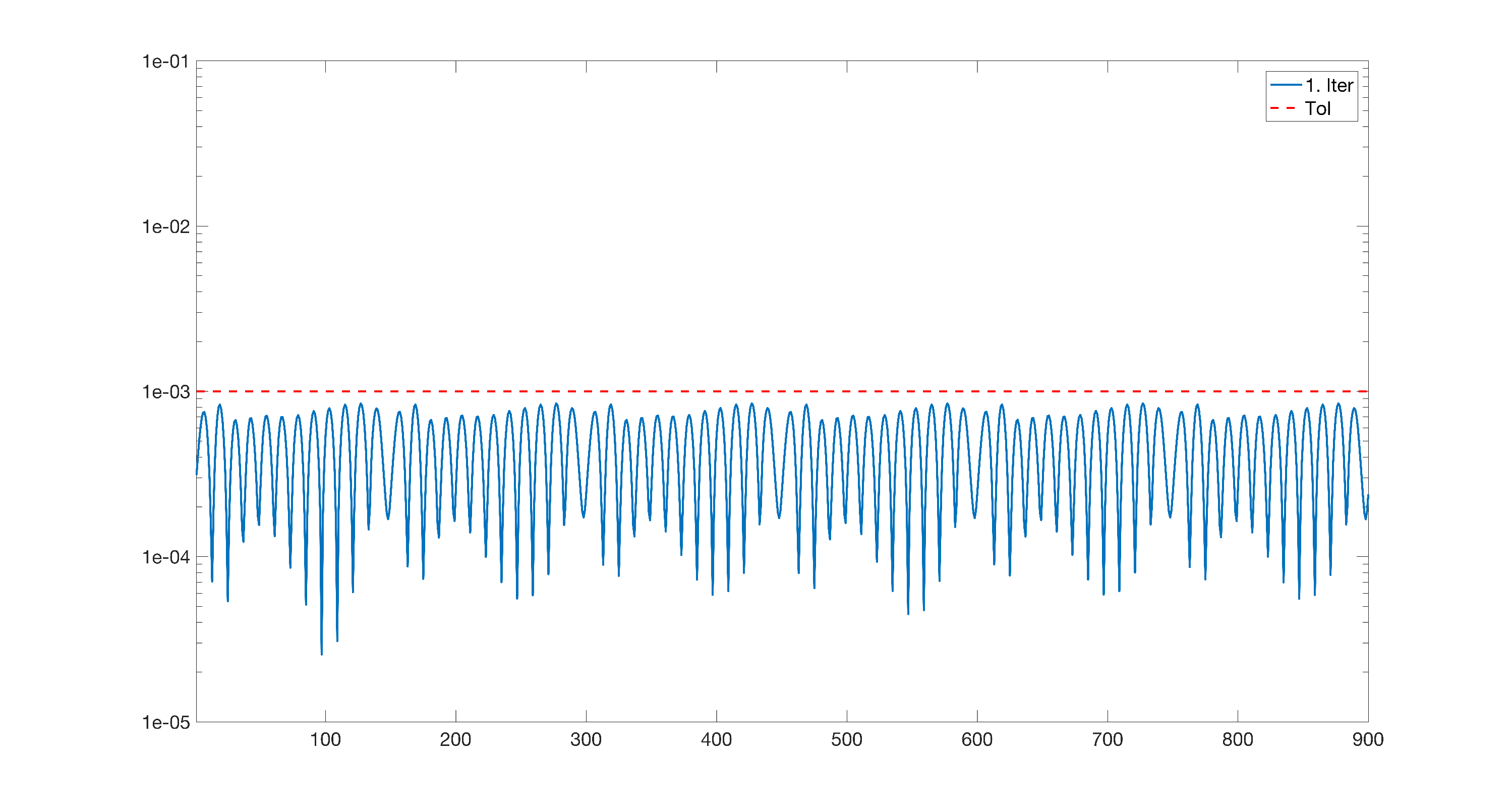}
        & \includegraphics[width=58mm]{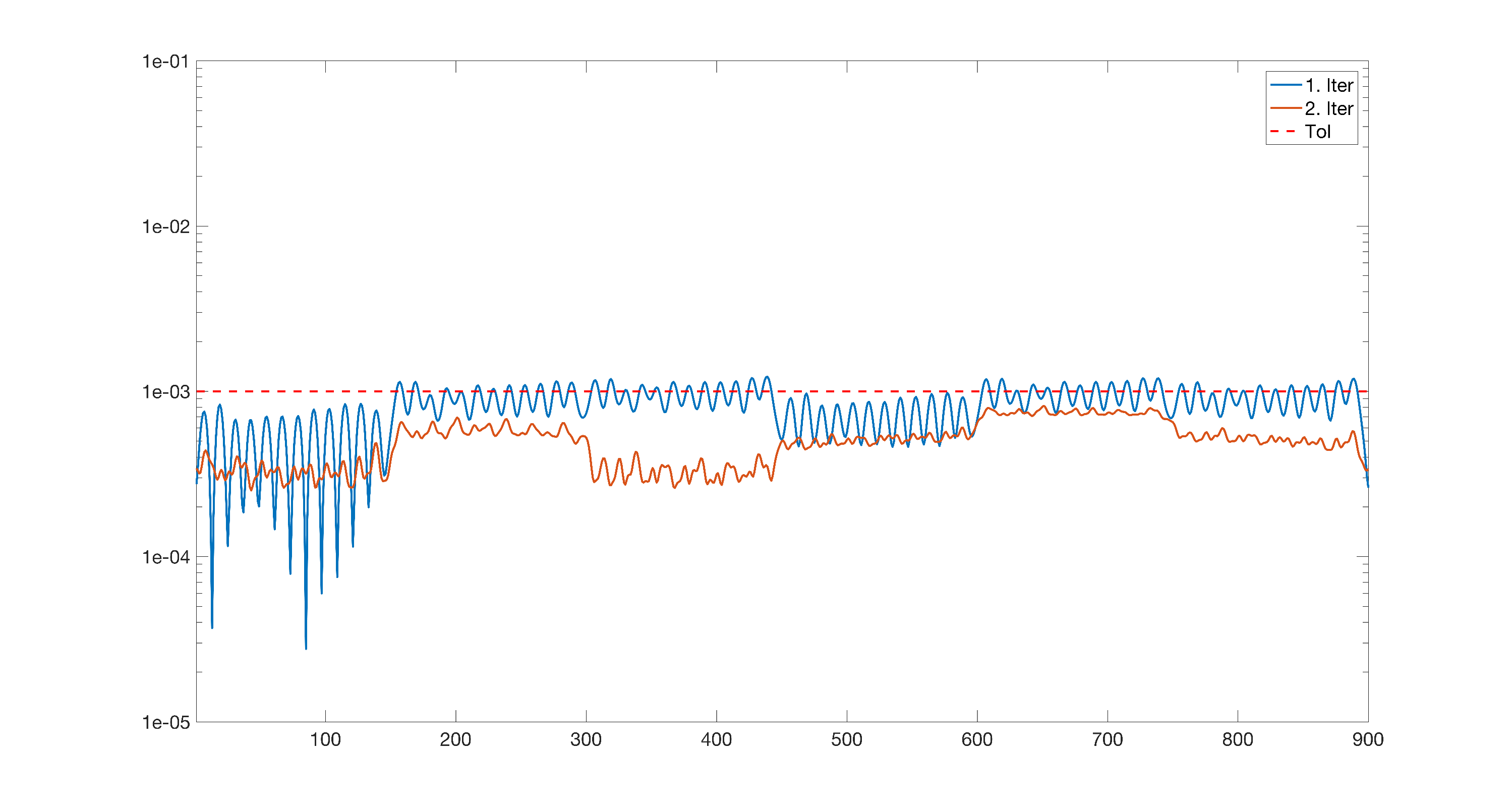}
        \\
          (a) \emph{sym}
        & (b) \emph{rot}
        \\[6pt]
          \includegraphics[width=58mm]{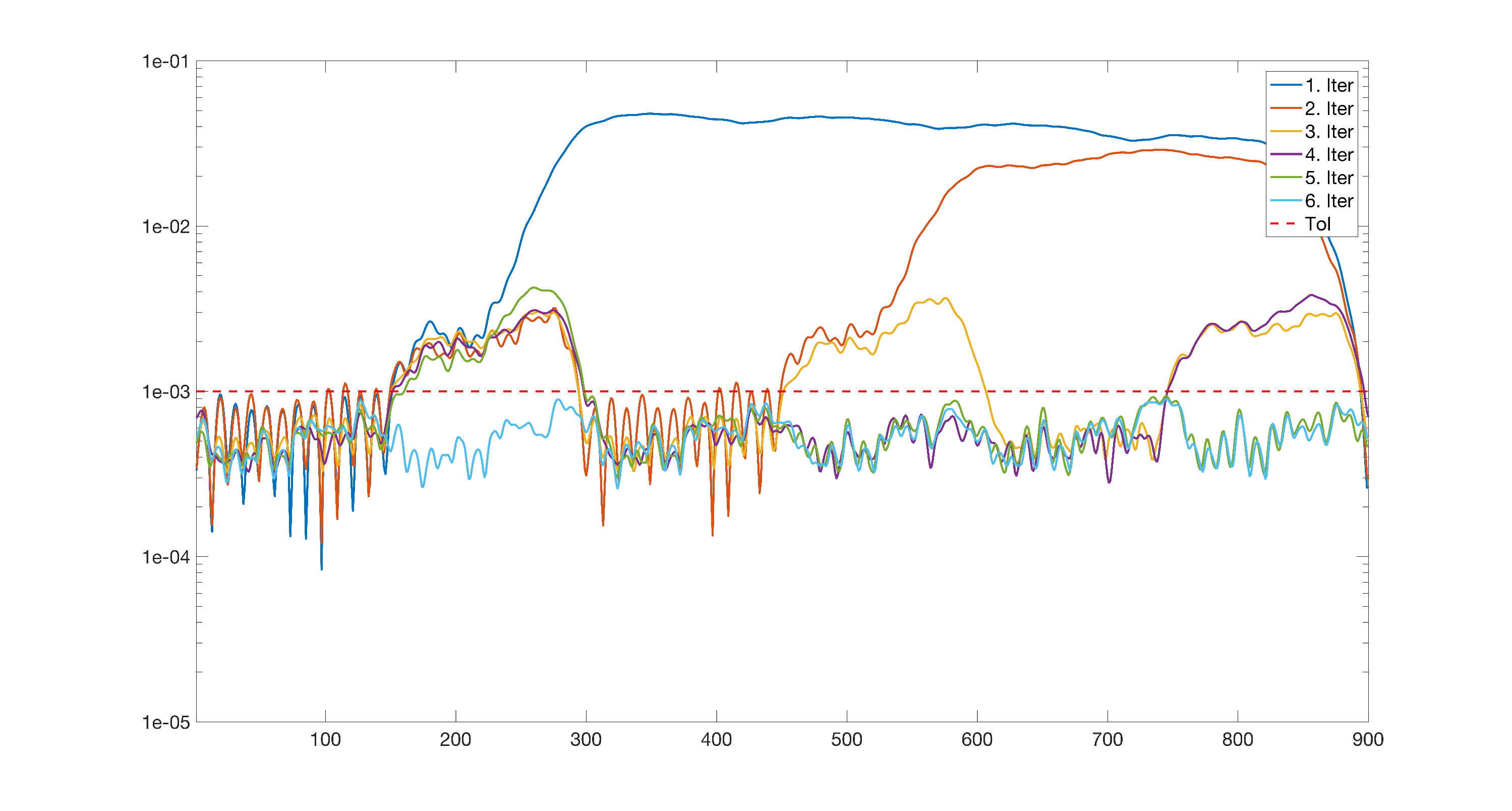}
        & \includegraphics[width=58mm]{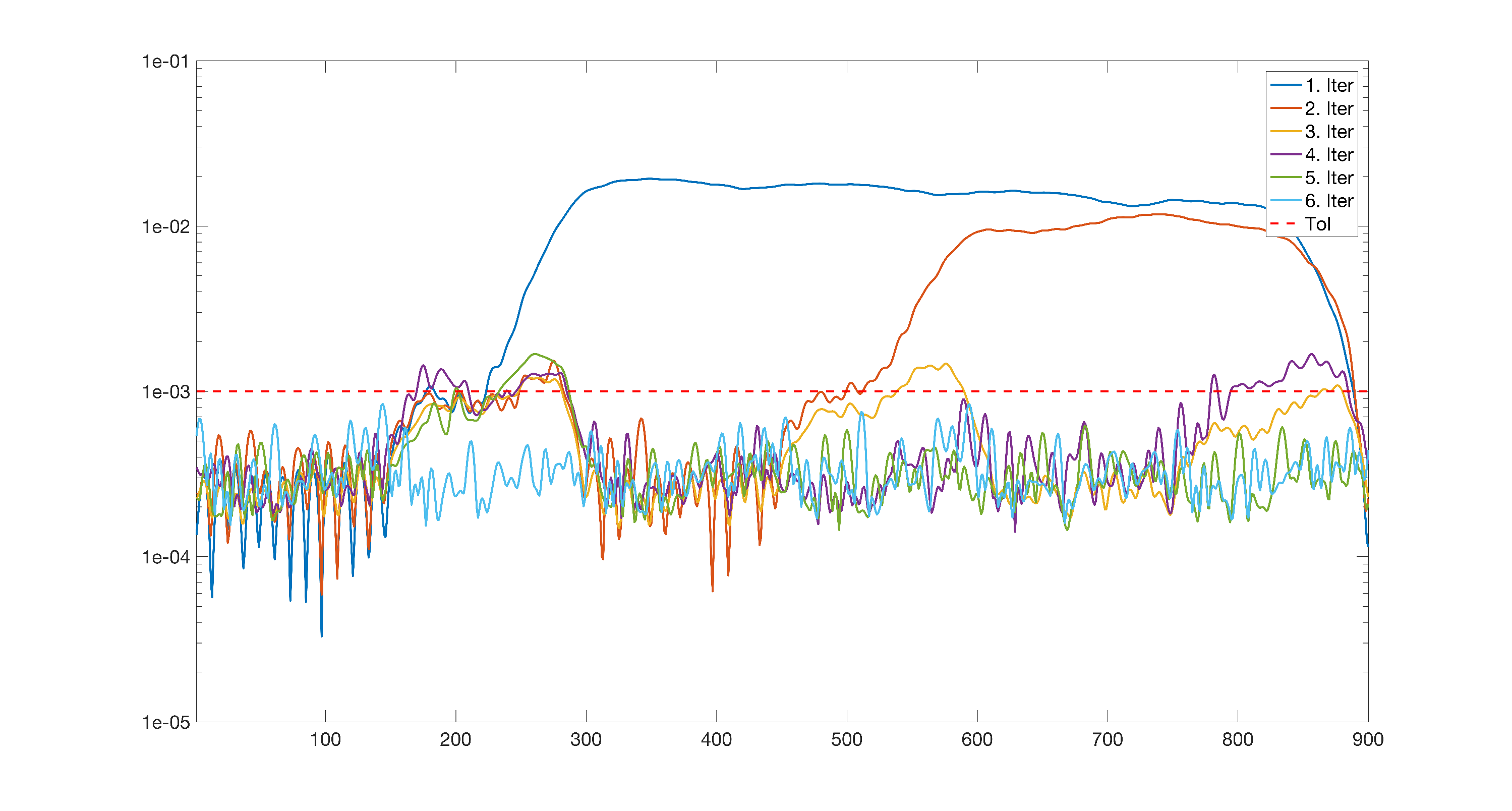}
        \\
          (c) \emph{stat}
        & (d) \emph{rot\_stat}
    \end{tabular}
    \caption{Error estimator in each iteration of the adaptive strategy using the sets $\mathcal K$.}
    \label{fig:res_K}
\end{figure}

In the second setting distributed index sets are considered. This allows for a
broader information collection since the snapshots are immediately distributed
over the whole range of $\vartheta$. We define the corresponding sets by
\begin{equation*} \mathcal M_{i} = [i, i+72, i+144, \ldots ],\quad i =
1,\ldots,72. \end{equation*} With these particular choices we get in both
settings $72$ sets with each $12$ indices and hence we can perform a fair
comparison of the two approaches. In our experiment it turned out that these
settings are a good trade between performance and accuracy. As the tolerance
for our adaptive strategy we use $\varepsilon = 10^{-3}$ which is sufficient
for most applications. The error estimator overestimated the real error by at
most one order of magnitude. This is a good result since it guarantees that the
model is not refined unnecessarily often.

\begin{figure}[H]
    \begin{tabular}{cc}
           \includegraphics[width=58mm]{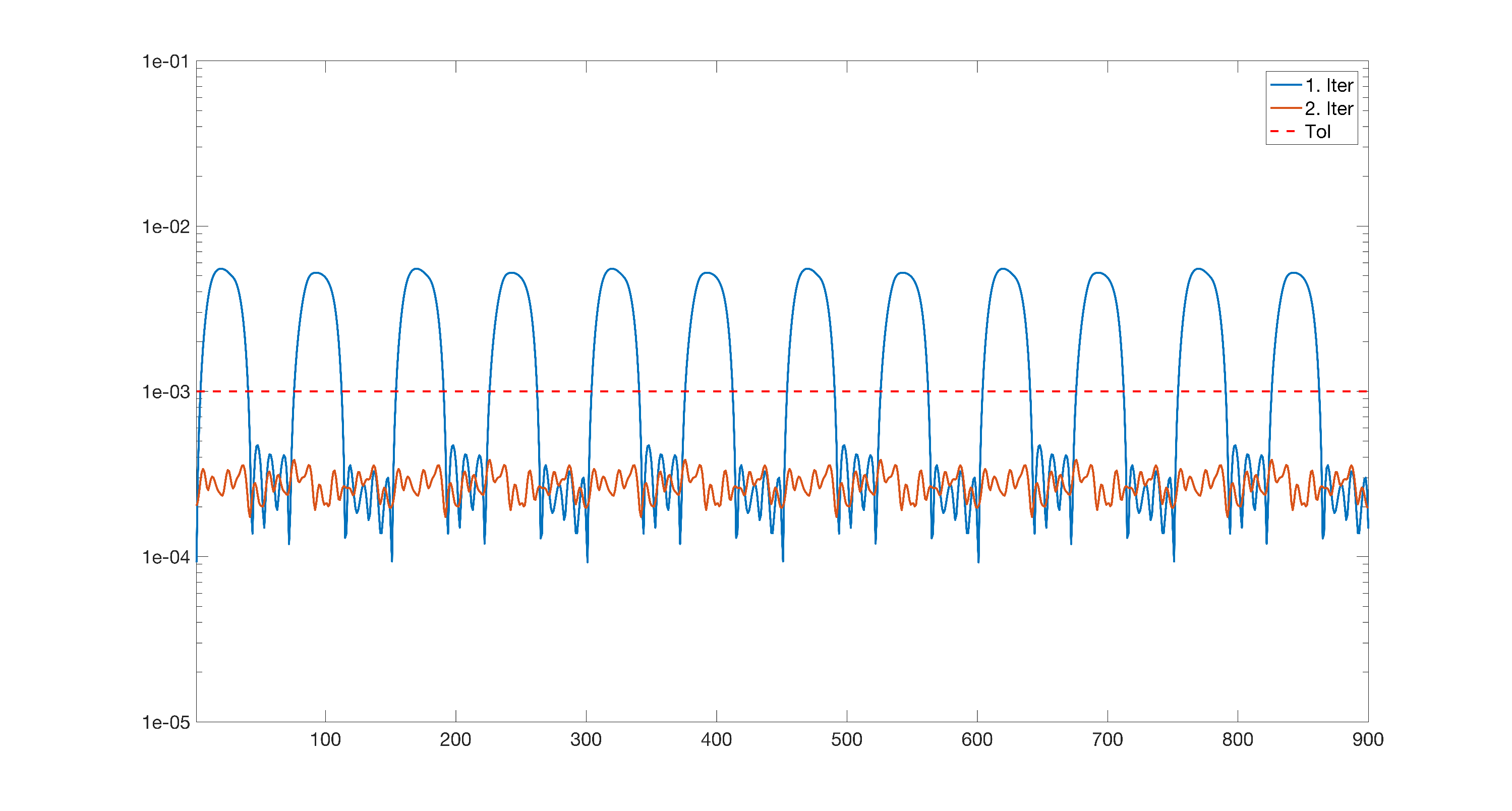}
        & \includegraphics[width=58mm]{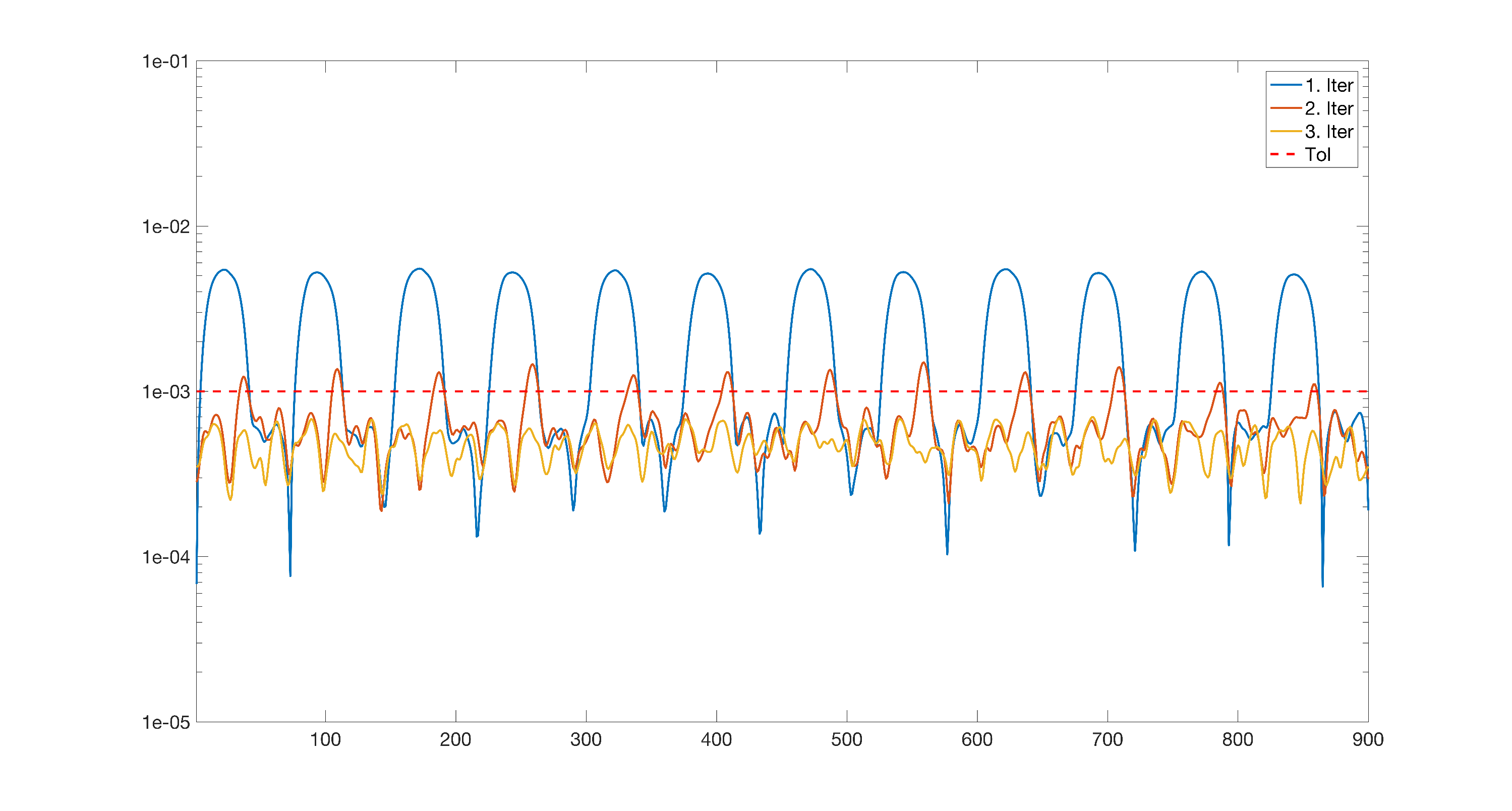}
        \\
          (a) \emph{sym}
        & (b) \emph{rot}
        \\[6pt]
          \includegraphics[width=58mm]{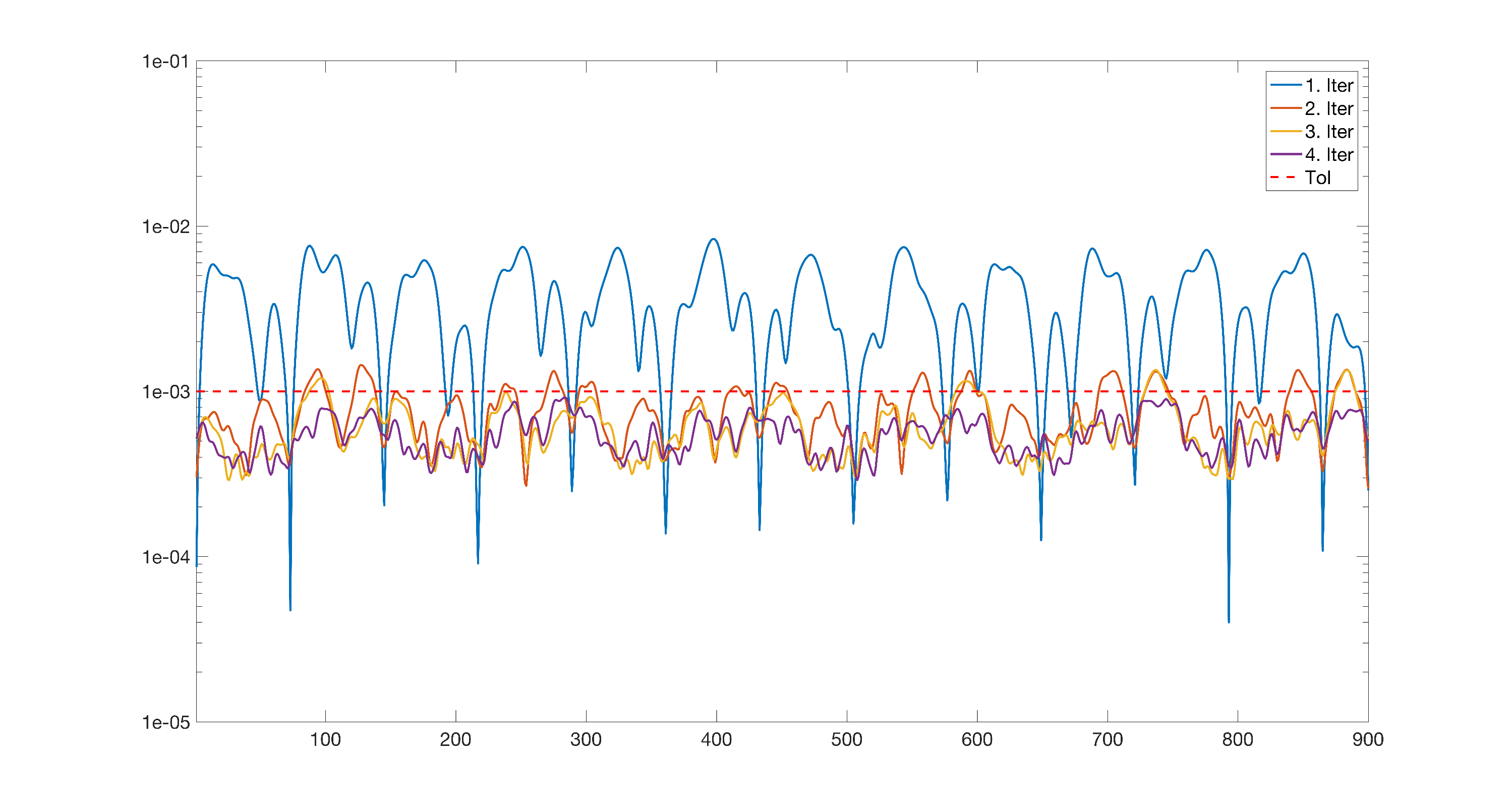}
        & \includegraphics[width=58mm]{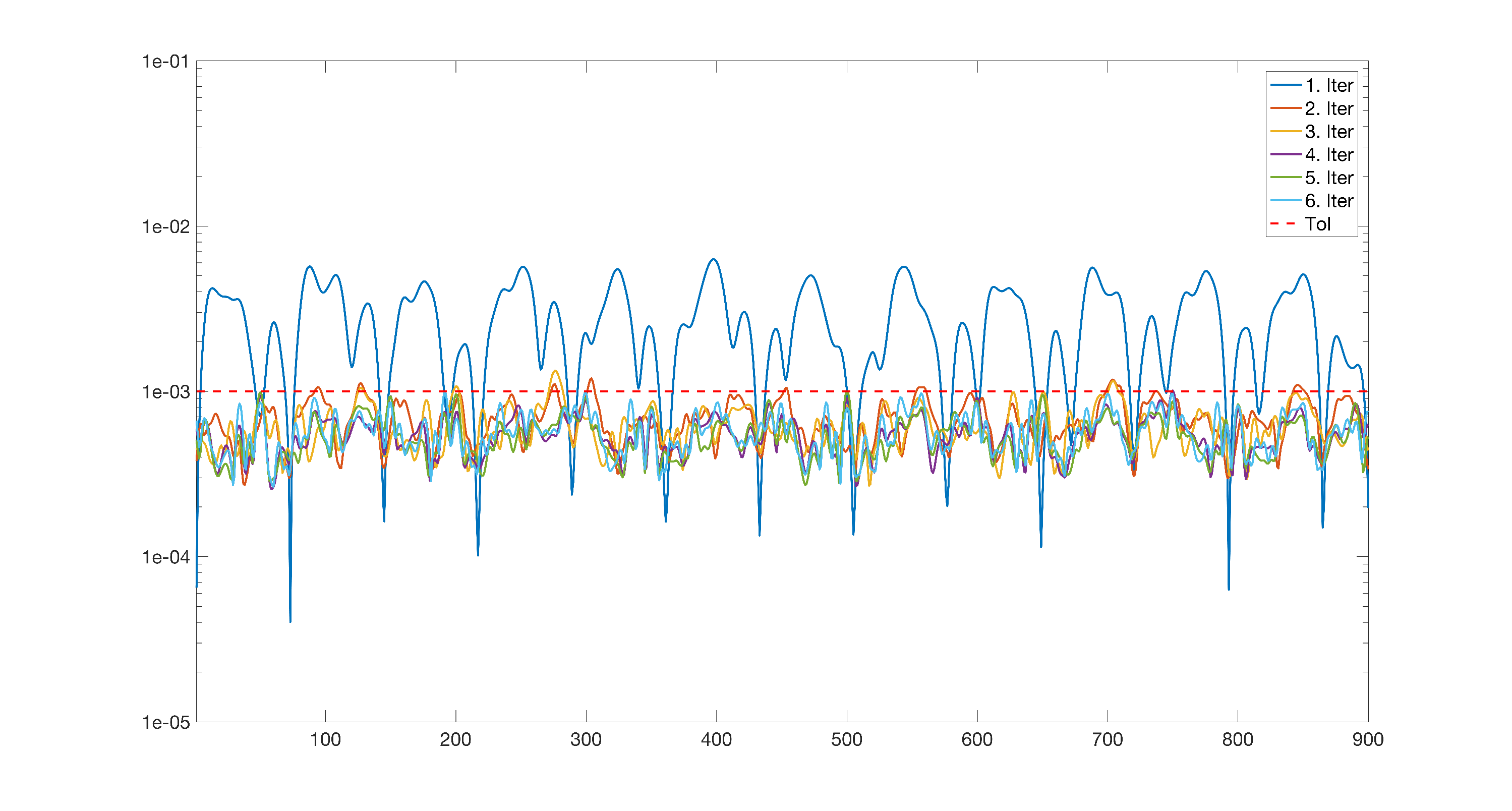}
        \\
          (c) \emph{stat}
        & (d) \emph{rot\_stat}
    \end{tabular}
    \caption{Error estimator in each iteration of the adaptive strategy using the sets $\mathcal M$.}
    \label{fig:res_M}
\end{figure}

The results for the two different index sets $\mathcal K$ and $\mathcal M$ are shown in Figure~\ref{fig:res_K}-\ref{fig:res_M}.
We plot the relative error $\Delta^{\mathrm{rel}}_{\mathbf a}$ of the $900$ angular position for every iteration of the adaptive algorithm. 
The actual error is omitted in the plots for visual clarity. 
By a dashed line the tolerance is indicated as a visual aide.
It can be seen how each of the two index sets have very different behaviors.

The set $\mathcal K$, which uses a kind of local refinement of the snapshots,
performs well for the \emph{sym} and \emph{rot}. In the behaviour of the error
it can be clearly seen in which region new snapshots were added. The error of
the region drops significantly. As was expected, the symmetric case only
requires one iteration since already a few snapshots contain enough information
to compute the full rotation. For the perturbation in the rotor a second set of
snapshots is required to push the error below the tolerance $\varepsilon$. A
clear accuracy difference between the poles can be observed which is reflected
by steps in the error. The settings \emph{stat} and \emph{rot\_stat} are not
handled too well by the local nature of the set $\mathcal K$. For each pole a
snapshot set is selected resulting in six iterations for both settings.

On the other hand the index set $\mathcal M$ shows the advantages of the global
nature of the snapshot selection. While for \emph{sym} and \emph{rot} this
results in more iterations as for the set $\mathcal K$ in the other two
settings, the benefits can be observed. In particular for \emph{stat}, a much
faster convergence of the adaptive algorithm is obtained. The error also has a
different nature. In each iteration the error is reduced uniformly over the
whole rotation and not only locally.

Next we have a look at the performance of the adaptive algorithm. For this we
compare computational time (wall clock time) of the different approaches. The
computational time is determined by the average over $10$ runs to flatten out
irregularities in the numerical realization. Note that the mesh and the
constant finite element matrices are precomputed, and hence are not reflected
in the run time. This does not influence the result since it is required by
both approaches and gives a more clear indication of the computational costs,
which is the main focus. In Table~\ref{tab:res_K}-\ref{tab:res_M} the results
are summarized.

\begin{table}
\caption{Performance summary for the different settings and index set $\mathcal K$.}
\label{tab:res_K} 
\begin{tabular}{p{1.5cm}p{1.75cm}p{1.75cm}p{2cm}p{1cm}p{1.5cm}p{1.5cm}}
\hline\noalign{\smallskip}
Setting          & FEM (sec.) & ROM (sec.) & Basis $(\statdom,\rotdom)$ & Iter. & Speedup & Overhead\\
\noalign{\smallskip}\svhline\noalign{\smallskip}
\emph{sym}       & $328.56$ & $\phantom{0}7.04$ & $(12,\phantom{0}6)$          & $1$   & $46.67$ & $37$\% \\
\emph{rot}       & $337.17$ & $12.76$ & $(12,\phantom{0}6)$          & $2$   & $26.42$ & $29$\%\\
\emph{stat}      & $341.28$ & $46.98$ & $(13,16)$         & $6$   & $\phantom{0}7.26$ & $41$\%\\
\emph{rot\_stat} & $328.30$ & $49.85$ & $(13,16)$         & $6$   & $\phantom{0}6.58$ & $47$\%\\
\noalign{\smallskip}\hline\noalign{\smallskip}
\end{tabular}
\end{table}

As already observed in the figures the two index sets have different strength
and weaknesses. When looking at the raw performance the speedup can vary
significantly. We get a factor of $46$ for \emph{sym} and go as low as $6$ for
\emph{rot\_stat} when using $\mathcal K$. Overall the speedup obtained by the
index set $\mathcal M$ is better since it does not exhibit too strong
variations but for \emph{sym} and \emph{rot} the local nature of $\mathcal K$
is significantly better.

What can be observed is that although more and more snapshots are being added
the dimension of the POD basis is almost the same for the different settings.
For all settings, only $13$ basis vectors for the stator are required, while
the number of basis vectors for the rotor depends on the problem setting. As
observed during the investigation of the decay of the eigenvalues, the settings
\emph{stat} and \emph{rot\_stat} require more basis vectors.

\begin{table}
\caption{Performance summary for the different settings and index set $\mathcal M$.}
\label{tab:res_M}
\begin{tabular}{p{1.5cm}p{1.75cm}p{1.75cm}p{2cm}p{1cm}p{1.5cm}p{1.5cm}}
\hline\noalign{\smallskip}
Setting          & FEM (sec.) & ROM (sec.) & Basis $(\statdom,\rotdom)$ & Iter. & Speedup & Overhead\\
\noalign{\smallskip}\svhline\noalign{\smallskip}
\emph{sym}       & $328.56$ & $14.36$ & $(13,\phantom{0}6)$        &  $2$    & $22.88$ & $38$\% \\
\emph{rot}       & $337.17$ & $20.60$ & $(13,\phantom{0}6)$        &  $3$    & $16.36$ & $34$\% \\
\emph{stat}      & $341.28$ & $32.74$ & $(13,16)$        & $4$    & $10.42$ & $44$\% \\
\emph{rot\_stat} & $328.30$ & $42.17$ & $(13,16)$        & $6$    & $\phantom{0}7.78$ & $37$\% \\
\noalign{\smallskip}\hline\noalign{\smallskip}
\end{tabular}
\end{table}

Lastly, we have a look at the distribution of the computational time. For the
FEM we require $300$-$350$ seconds to complete the simulation. Considering,
that $900$ linear systems have been solved, the average time for solving one
linear system is $0.33$-$0.38$ second. Multiplying this with the number of
computed snapshots we can determine the overhead in the computation introduced
by the basis computation, ROM generation and evaluation of the error estimator.
It turns out that the overhead is between $30$\% and $50$\%. Here, efficient
methods can be investigated that can update existing POD basis more
efficiently. Further, we utilized the SVD method to compute the basis vectors.
This was done for stability reasons so that during the tests we do not
encounter problems. Alternatively, the formulation using the eigenvalue
representation can be investigated. This may give a boost but the overall
performance will lead to similar conclusions.

Overall the numerical results are very pleasing and the model order reduction
is very effective. The speedup might not be as large when compared to
strategies utilizing an online-offline decomposition. On the other hand the
presented framework can be used to directly replace the simulation routine and
no further adjustment need to be done to existing code. Especially, when
embedding into optimization solvers or sampling methods like Monte Carlo
simulations this can turn out as a big benefit. Further, when the ROM are
utilized a significant reduction in memory usage can be achieved. This then
allows large simulation even on moderate hardware like desktop PCs.

\section{Conclusion}
\label{sec:con}

We developed an adaptive POD snapshot sampling strategy targeted at model order
reduction of electric rotating machines. A detailed description of the required
component is provided. In the numerical results different strategies of
snapshot sampling were investigated and compared. The method proved to be very
efficient in reducing the computational cost of symmetric and non-symmetric
machines.

\begin{acknowledgement}
This work is supported by the German BMBF in the context of the SIMUROM project (grant number 05M2013), by the ‘Excellence Initiative’ of the German Federal and State Governments and the Graduate School of Computational Engineering at TU Darmstadt.
\end{acknowledgement}

\end{document}